\documentclass[11pt]{article}
\usepackage{authblk}

\usepackage[ansinew]{inputenc}
\usepackage{amsmath,amssymb,amsthm}
\usepackage{subfig}
\usepackage{graphicx}
\usepackage{multirow}
\usepackage{color}
\usepackage[bbgreekl]{mathbbol}

\usepackage[]{algorithm2e}

\DeclareMathAlphabet{\mathpzc}{OT1}{pzc}{m}{it}

\newcommand{\coeff}[1]{\ensuremath{\lowercase{#1}}}

\renewcommand{\vec}[1]{\ensuremath{\mathbf{#1}}}
\newcommand{\dist}{\ensuremath{\mathpzc{f}}}
\newcommand{\ham}{\ensuremath{\mathcal{H}}}

\DeclareMathOperator{\curl}{curl}

\newcommand{\du}{\,\mathrm{d}}
\newcommand{\R}{\mathbb{R}}

\newcommand{\im}{\mathrm{i}}
\newcommand{\e}{\mathrm{e}}

\newcommand{\xb}{\mathbf{x}}
\newcommand{\vb}{\mathbf{v}}
\newcommand{\ub}{\mathbf{u}}
\newcommand{\uu}{\mathbf{u}}
\newcommand{\eb}{\mathbf{e}}

\newcommand{\Eb}{\mathbf{E}}
\newcommand{\Bb}{\mathbf{B}}
\newcommand{\bb}{\mathbf{b}}

\newcommand{\jac}{\ensuremath{\mathcal{J}}}

\newcommand{\BB}{\ensuremath{\mathbb{B}}}
\newcommand{\MM}{\ensuremath{\mathbb{M}}}

\newcommand{\C}{\ensuremath{\mathbb{C}}}
\newcommand{\D}{\ensuremath{\mathbb{D}}}
\newcommand{\G}{\ensuremath{\mathbb{G}}}

\newcommand{\DD}{\ensuremath{\mathsf{D}}}
\newcommand{\f}{\ensuremath{\dist}}

\newcommand{\E}{\ensuremath{\vec{E}}}
\newcommand{\J}{\ensuremath{\vec{J}}}
\newcommand{\X}{\ensuremath{\vec{X}}}
\newcommand{\V}{\ensuremath{\vec{V}}}
\newcommand{\U}{\ensuremath{\vec{U}}}

\newcommand{\Lab}{\ensuremath{\boldsymbol{\Lambda}}}
\newcommand{\LaB}{\ensuremath{\mathbb{\Lambda}}}

\newcommand{\cb}{\ensuremath{\vec{\coeff{B}}}}

\newcommand{\ce}{\ensuremath{\vec{\coeff{E}}}}
\newcommand{\cf}{\ensuremath{\vec{\coeff{F}}}}
\newcommand{\cj}{\ensuremath{\vec{\coeff{J}}}}

\newcommand{\crho}{\ensuremath{\boldsymbol{\varrho}}}

\newcommand{\dd}{\ensuremath{\mathrm{d}}}

\newcommand{\ii}{\text{iter}}
\newcommand{\me}{\tilde{e}}
\newcommand{\mb}{\tilde{b}}

\title{Energy-conserving time propagation for a geometric particle-in-cell Vlasov--Maxwell solver}

\author{Katharina Kormann, Eric Sonnendr\"ucker}
\affil{Max Planck institute for Plasma Physics, Boltzmannstr.~2, 85748 Garching, Germany and Technical University of Munich, Department of Mathematics, Boltzmannstr.~3, 85748 Garching, Germany.}
\date{}
\begin{document}

\maketitle




\begin{abstract}
This paper discusses energy-conserving time-discretizations for finite element particle-in-cell discretizations of the Vlasov--Maxwell system. A geometric spatially discrete system can be obtained using a standard particle-in-cell discretization of the particle distribution and compatible finite element spaces for the fields to discretize the Poisson bracket of the Vlasov--Maxwell model (see Kraus et al., J Plasma Phys 83, 2017). In this paper, we derive energy-conserving time-discretizations based on the discrete gradient method applied to an antisymmetric splitting of the Poisson matrix. Firstly, we propose a semi-implicit method based on the average-vector-field discretization of the subsystems. Moreover, we devise an alternative discrete gradient that yields a time discretization that can additionally conserve Gauss' law. Finally, we explain how substepping for fast species dynamics can be incorporated.
\end{abstract}



\section{Introduction}

Particle-in-cell simulations are widely used in the plasma community to solve the Vlasov--Maxwell's equations due to their ease of implementation and their favorable scaling properties in higher dimensions.
Recently, a systematic derivation of geometric particle-in-cell methods has been proposed by Kraus, Kormann, Morrison, \& Sonnendr\"ucker \cite{Kraus17}. The derivation is based on compatible finite elements for the fields and a standard particle-in-cell ansatz for the particles. The derived semi-discrete system conserves Casimir invariants of the system such as discrete versions of $\nabla \cdot \Bb = 0$ and $\nabla \cdot \Eb = \rho$. For the time discretization, a Hamiltonian splitting method was proposed which yields an explicit scheme that conserves Gauss' law over time, however, only a modified energy. While one step of this method is very efficient, the time step is restricted by stability constraints (cf.~the Appendix). Moreover, the splitting only yields an explicit scheme as long as the coordinate system is orthogonal that is it does not yield a simple scheme for the case of curvilinear coordinates. This motivates our investigations on alternative time-stepping schemes for the geometric electromagnetic particle-in-cell method. 

In this paper, we devise alternative temporal discretizations based on discrete gradients. The discrete gradient method is a general framework to design energy-conserving time discretizations for conservative partial differential equations in skew-symmetric form and was first introduced by McLachlan, Quispel, \& Robidoux \cite{McLachlan99}. Several special cases have been devised, in particular the average-vector-field method introduced by Celledoni et al.~\cite{Celledoni12}. Applying discrete gradients to the full Vlasov--Maxwell system results in a heavily nonlinear scheme. On the other hand, the method of discrete gradients can be applied after a splitting of the equations that respects the skew-symmetry without loss of energy conservation. As a first scheme, we propose a semi-implicit method that applies the average-vector field method to the subsystems that cannot be solved analytically. The scheme is only implicit in the field solver, whence the computational overhead is relatively small compared to the explicit method. On the other hand, this method does not conserve Gauss' law. We therefore devise a second scheme where we reduce the splitting and derive an alternative discrete gradient that conserves Gauss' law in addition to the energy.  The resulting system nonlinearly couples the particle and field equations and therefore needs to be solved in a nonlinear iteration. In a simulation with fast electrons and slower ions, substepping for the trajectories of the faster species can be crucial to reflect the multiscale nature of the system. Such a substepping technique can be incoorporated in our implicit scheme.

In the plasma physics community, energy-conserving particle-in-cell methods have also been developed, mostly with a finite difference description of the fields. Markidis \& Lapenta \cite{Markidis11} have devised the so-called EC-PIC method for the Vlasov--Maxwell system that is fully nonlinear: The method uses a finite difference description of the fields on a Yee grid and employs differencing by the implicit midpoints rule both in space and time. With some rearrangements of the equations they yield an implicit formulation for the update of the fields and an average velocity. The method conserves energy but Gauss' law is not preserved over time. Later Lapenta derived a semi-implicit version, the so-called energy-conserving semi-implicit particel-in-cell method (ECSIM) \cite{Lapenta17}, that only requires an implicit field solver. When transferring the ECSIM method to a finite element formulation of the fields, the resulting method is very close to our energy-conserving semi-implicit scheme derived from the average-vector-field method. This has been investigated in \cite{Perse17}. Compared to the ECSIM method, we follow a more systematic derivation that ensures second-order accuracy which is lost in one part of the ECSIM method as pointed out in \cite{Perse17}. Similar to the average-vector-field method derived in this work, the ECSIM method is energy-conserving but does not satisfy Gauss' law. On the other hand, Gauss' law can be reinforced by different Lagrange multiplier techniques as proposed by Marder \cite{Marder87}, Langdon \cite{Langdon92} and by Munz, Omnes, Schneider, Sonnendr\"ucker, \& Voss \cite{Munz.Omnes.Schneider.Sonnendrucker.Voss.1999.jcp}. However, these techniques are not compatible with the energy conservation in our scheme. Recently, Chen \& Toth \cite{Chen19} have therefore proposed a different procedure for the ECSIM algorithm that uses a correction of the particle positions instead and respects the energy conservation. However, this method is build on linearized shape functions and is therefore only a first order correction.

Chen, Chac\'on, \& Barnes have developped fully implicit particle-in-cell methods that conserve both energy and Gauss' law for the electrostatic Vlasov--Amp\`ere model \cite{Chen11} as well as the reduced electromagnetic Vlasov--Darwin model \cite{Chen15}.

The outline of the paper is as follows: In the next section, we introduce the Vlasov--Maxwell model and the geometric electromagnetic particle-in-cell (GEMPIC) framework for its spatial semi-discretization. Section \ref{sec:disgrad_general} introduces the discrete gradient method and a splitting of the Poisson matrix of the semi-discrete Vlasov--Maxwell system and Section \ref{sec:avf} devises an energy-conserving average-vector-field method for the split equations. In the subsequent section \ref{sec:disgrad}, we explain how the scheme can be modified to conserve as well Gauss' law. Numerical experiments on simple test problems presented in Section \ref{sec:numerics} confirm the conservation properties of the new methods.

\section{Geometric electromagnetic particle in cell}

In the GEMPIC framework \cite{Kraus17}, the Vlasov--Maxwell equations are discretized by a standard particle-in-cell ansatz for the distribution function and compatible finite elements for the fields. The spatial semi-discretization is derived from a semi-discretization of the Hamiltonian and the Poisson bracket. In the following, we revise this semi-discretization as a starting point for the time-discretizations proposed in this paper.

\subsection{The Vlasov--Maxwell system}

A kinetic description of a plasma models a species $s$ of particles with charge $q_s$ and mass $m_s$ by a distribution function $\f_s$ in phase-space that evolves according to the Vlasov equation
\begin{equation*}
\frac{\partial \f_s}{\partial t} + \vb \cdot \nabla_{\xb} \f_s + \frac{q_s}{m_s} \left( \Eb + \vb \times \Bb \right) \cdot \nabla_{\vb} \f_s = 0,
\end{equation*} 
where $\Eb$ and $\Bb$ denote the external and self-consistent electric and magnetic fields. The advection equation is coupled to the Maxwell's equations for the self-consistent fields
\begin{subequations}\label{eq:vlasov_maxwell_equations}
\begin{align}
&\frac{\partial \Eb}{\partial t} - \nabla \times  \Bb = - \J, \quad 
\label{eq:ampere} \\
&\frac{\partial \Bb}{\partial t}  + \nabla  \times \Eb = 0, \label{eq:faraday}
\\
&\nabla \cdot \Eb = \rho, \quad \label{eq:gauss}\\
&\nabla \cdot \Bb = 0, \label{eq:grad_b}
\end{align}
\end{subequations}
where the charge density $\rho$ and the current density $\J$ are defined as velocity moments of the distribution functions
\begin{equation*}
\rho = \sum_s q_s \int \f_s \du \vb, \text{   and   } \J = \sum_{s} q_s \int \vb \f_s \du \vb.
\end{equation*}
As for any hyperbolic conservation law, the solution stays constant along the characteristic equations, which are defined as the following system of ordinary differential equations,
\begin{equation}\label{eq:characteristics}
\frac{\du \xb}{\du t} = \vb, \quad \frac{\du \vb}{\du t} = \frac{q_s}{m_s} \left( \Eb(\xb,t) + \vb \times \Bb(\xb,t)\right).
\end{equation}
The following Hamiltonian defines the total energy of the system
\begin{align}
\ham = \sum_{s} \frac{m_s}{2} \int |{\vb}|^{2} \, \f_s (t,\xb,\vb) \du \xb \du\vb + \frac{1}{2} \int \Big( |{\Eb (t, \xb)}|^{2} + |{\Bb (t, \xb)}|^{2} \Big)  \du \xb .
\label{eq:hamiltonian_vlasov_maxwell}
\end{align}


\subsection{Compatible finite element discretization of the field equations}

\begin{figure}
\centering
\includegraphics[scale=1.0]{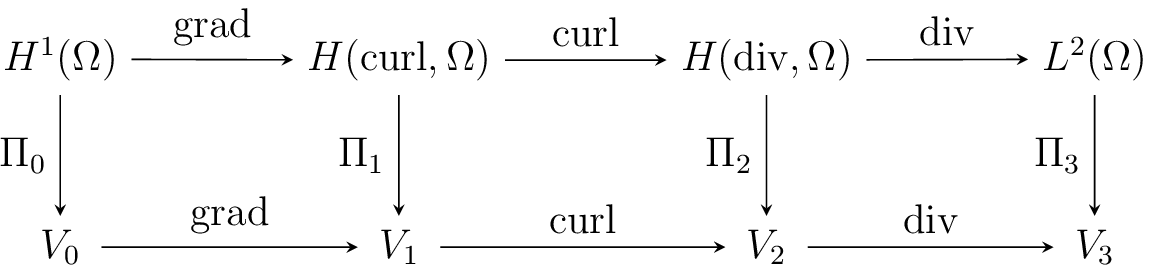}
\caption{Discrete de Rham complex for the spaces of electromagnetics.}\label{fig:de_rham}
\end{figure}

The Maxwell's equations themselves posses a rich structure of conservation properties. Especially, the spaces of electromagnetics form a de Rham complex (cf.~the first line in Figure \ref{fig:de_rham}) with $\Eb \in H(\curl,  \Omega)$, $\Bb, \J \in H(\text{div},\Omega)$, and $\rho \in L^2(\Omega)$. The theory of finite element exterior calculus \cite{Arnold.Falk.Winther.2006.anum,Arnold.Falk.Winther.2010.bams} explains how these properties can be preserved in a finite element discretization: The discrete spaces are chosen in such a way that they form the commuting diagram shown in Figure \ref{fig:de_rham} with the continuous spaces.   We use $\Eb_h \in V_1$ and $\Bb_h \in V_2$. Denoting the $V_1$-basis functions by $ \Lab^{1}$ and the $V_2$ basis by $ \Lab^{2}$, we get the following discrete representation of the fields 
\begin{eqnarray}\label{eq:discrete_eb}
\Eb_h(\xb) = \sum_{i=1}^{3N_1} e_i(t)  \Lab^{1}_{i} (\xb), \quad
\Bb_h(\xb) = \sum_{i=1}^{3N_2} b_i(t) \Lab^{2}_{i} (\xb),
\end{eqnarray}
where $\eb = ( e_i)$ and $\bb = (b_i)$ are the degrees of freedom in the semi-discretization. A compatible finite element discretization of the Maxwell's equation can be obtained based on the ansatz  \eqref{eq:discrete_eb} treating \eqref{eq:ampere} and \eqref{eq:gauss} in weak and \eqref{eq:faraday} and \eqref{eq:grad_b} in strong form,
\begin{subequations}
\begin{align*}
\MM_{1} \frac{\dd \ce}{\dd t} - \C^\top \MM_{2} \cb &= -\cj , \\
\frac{\dd \cb}{\dd t} +\C \ce &= 0 , \\
\G^\top \MM_{1} \ce &= \crho  , \\
\D \cb &= 0  ,
\end{align*}
\end{subequations}
where $\MM_{1/2}$ are the finite element mass matrices for the basis functions $\Lab^{1/2}$, respectively, with elements $(\MM_{1/2})_{k\ell} = \int \Lab^{1/2}_{k}(\xb)\Lab^{1/2}_{\ell}(\xb) \du \xb$, and $\D$, $\G$, and $\C$ represent the discrete divergence, gradient, and curl operators which satisfy $\text{Im} \G = \text{Ker} \C$ and $\text{Im} \C = \text{Ker}  \D$ as their continuous counterparts.
There are various types of compatible finite element discretizations. In our work, we use spline finite elements of various order as proposed by Buffa, Rivas, Sangalli, \& V\'asquez \cite{Buffa:2011}. 
The charge and current density are tested with the corresponding basis functions to obtain the vectors $\crho$ and $\cj$ for the right-hand-side.

\subsection{Coupling to the particles}

A particle discretization represents the distribution function $\f_s$ by a number $N_{p_s}$ of particles which in turn are represented by the dynamic variables $(\xb_a, \vb_a)$, its coordinates in phase-space, as well as a weight $w_a$ which is fixed over time in our formulation. The particle distribution is reconstructed by the Klimontovich distribution
\begin{equation}
\f_{s,h} (\xb,\vb,t) = \sum_{a=1}^{N_{p_s}} \, w_a \, \delta \big( \xb - \xb_a (t) \big) \, \delta \big( \vb - \vb_a (t) \big) .
\label{eq:vlasov_maxwell_distribution_discrete}
\end{equation}
This representation is suitable in a finite-element discretization where the velocity moments are only needed in the weak form. In the resulting semi-discretized system, the dynamic variables $\ub^\top = \left( \X^\top, \, \V^\top, \, \eb^\top, \, \bb^\top\right)$ are given by the phase space positions of the particles (for all species), $(\X^\top, \, \V^\top)$, 
and the degrees of freedom of the fields, $\eb = ( e_i)$ and $\bb = (b_i)$. The integrals over the test functions times the charge and current density, that are needed for the right-hand-side of the Amp\`ere's law and the electric Gauss' law, can be computed as follows
\begin{equation*}
\rho_i = \int \Lab_i^{1}(\xb) \rho(\xb) \du \xb = \sum_s q_s \int \Lab_i^{1}(\xb) \f_{s,h} (\xb,\vb,t) \du \xb \du \vb = \sum_{s} q_s \sum_{N_{p_s}} w_a   \Lab_i^{1}(\xb_a)
\end{equation*}
and analogously
\begin{equation*}
\cj_i = \int  \Lab_i^{2}(\xb) \J(\xb) \du \xb = \sum_s q_s \sum_{N_{p_s}} \vb_{a} w_a   \Lab_i^{2}(\xb_a).
\end{equation*}
The phase-space coordinates evolve in time according to the characteristic equations \eqref{eq:characteristics} which are semi-discretized as
\begin{equation*}
\frac{\du \xb_a}{\du t} = \vb_a, \quad \frac{\du \vb_a}{\du t} = \frac{q_a}{m_a} \left( \Eb_h(\xb_a) + \widehat{\Bb}_{h} (\xb_a, t) \vb_a \right),
\end{equation*}
where the matrix $\widehat{\Bb}_{h} (\xb_a, t)$ represents the crossproduct of a vector with $\Bb_{h}(\xb_a, t)$ such that \\ 
$\widehat{\Bb}_{h} (\xb_a)\vb_a=\vb_a \times \Bb_h(\xb_a)$ and can be computed as
\begin{align*}
\widehat{\Bb}_{h} (\xb_a, t) = \sum_{i=1}^{N_2} 
\begin{pmatrix}
\hphantom{-} 0 & \hphantom{-} b_{i,3}(t) \, \Lambda^{2,3}_{i} (\xb_a) & - b_{i,2}(t) \, \Lambda^{2,2}_{i} (\xb_a) \\
- b_{i,3}(t) \, \Lambda^{2,3}_{i} (\xb_a) & \hphantom{-} 0 & \hphantom{-} b_{i,1}(t) \, \Lambda^{2,1}_{i} (\xb_a) \\
\hphantom{-} b_{i,2}(t) \, \Lambda^{2,2}_{i} (\xb_a) & - b_{i,1}(t) \, \Lambda^{2,1}_{i} (\xb_a) & \hphantom{-} 0 \\
\end{pmatrix} .
\end{align*} 
In order to write the full semi-discrete system in matrix-vector form, we define the following matrices:  The diagonal matrices $\MM_q\in\R^{3N_p\times 3N_{p}}$ and $\MM_m\in \R^{3N_p\times 3N_{p}}$ collect the particle charges $q_s w_a$ or particles masses $m_s w_a$ on the diagonal and $\LaB^{1/2}(\X) \in \R^{3N_p\times 3N_{1/2}}$, the matrix containing the 
value of all basis functions on each particle position. Further, we denote by $\BB(\X,\cb) \in \R^{(3N_p)\times(3N_p)}$ the matrix that consists of one $3 \times 3$ block $\hat \Bb_h(x_a,t)$ for each particle. This yields the following set of semi-discrete equations of motion
\begin{subequations}\label{eq:vlasov_maxwell_equations_of_motion}
\begin{align}
\label{eq:vlasov_maxwell_equations_of_motion_x}
\dot{\X}
&= \V ,
\\
\label{eq:vlasov_maxwell_equations_of_motion_v}
\dot{\V}
&= \MM_p^{-1}\MM_q \big( \LaB^1 (\X) \ce + \BB (\X,\cb) \V \big) ,
\\
\label{eq:vlasov_maxwell_equations_of_motion_e}
\dot{\ce}
&= \MM_{1}^{-1} \big( \C^\top \MM_{2} \cb (t) - \LaB^1 (\X)^\top \MM_q \V \big) ,
\\
\label{eq:vlasov_maxwell_equations_of_motion_b}
\dot{\cb}
&= - \C \ce (t) .
\end{align}
\end{subequations}
Moreover, the semi-discretization of the Hamiltonian reads
\begin{align*}
\hat{\mathcal{H}}(\uu)
&= \tfrac{1}{2} \, \V^{\top}  \MM_{p} \V
 + \tfrac{1}{2} \, \ce^{\top} \MM_{1} \ce
 + \tfrac{1}{2} \, \cb^{\top} \MM_{2} \cb .
\end{align*}

In Section 4 of \cite{Kraus17}, these semi-discrete equations are derived from a semi-discretization of the Poisson bracket for the Vlasov--Maxwell system. This derivation leads to a representation of the equations of motion as a Hamiltonian system of the form
\begin{equation}\label{eq:eqns_of_motion_Poisson}
\dot{\uu} = \jac (\uu) \, \DD_{\uu} \hat \ham (\uu) ,
\end{equation}
where the skew-symmetric so-called Poisson matrix $\jac(\uu)$ is given as
\begin{equation}\label{eq:poisson_matrix}
\jac (\uu) =
\begin{pmatrix}
0 & \MM_p^{-1} & 0 & 0 \\
- \MM_p^{-1} & \MM_p^{-1} \MM_q \BB (\X,\cb) \, \MM_p^{-1} & \MM_p^{-1} \MM_q \LaB^1(\X) \MM_{1}^{-1} & 0 \\
0 & - \MM_{1}^{-1} \LaB^1 (\X)^\top \MM_q \MM_p^{-1} & 0 &  \MM_{1}^{-1} \C^\top\\
0 & 0 & - \C \MM_{2}^{-1} & 0 \\
\end{pmatrix} 
\end{equation}
The total derivative of the Hamiltonian can be computed to be $$\DD \hat \ham (\uu) = \left( 0 , \, \left(\MM_{p} \V\right)^\top , \, \left(\MM_{1} \ce\right)^\top , \, \left(\MM_{2} \cb\right)^\top\right)^{\top}.$$
In section \ref{sec:disgrad_general}, we will construct energy-conserving time discretizations based on this special form of the semi-discretization.

\subsection{The explicit Hamiltonian splitting}

Finally the system of equations \eqref{eq:vlasov_maxwell_equations_of_motion} has to be discretized in time. In \cite{Kraus17}, the temporal discretization is based on the form \eqref{eq:eqns_of_motion_Poisson} of the evolution equation combined with a splitting of the equations by splitting the Hamiltonian as 
\begin{equation*}
\hat \ham_{p_i} = \tfrac{1}{2}m_s \,  \sum_{a=1}^{N_p} w_p v_{i,a}^2, \, i=1,2,3,  \quad
\hat \ham_{\eb} = \tfrac{1}{2} \, \ce^{\top} \MM_{1} \ce, \quad 
\hat \ham_{\bb} = \tfrac{1}{2} \, \cb^{\top} \MM_{2} \cb.
\end{equation*}
This yields five sets of explicit equations (cf.~\cite[Sec.~5.1]{Kraus17}). The discrete system then still conserves Gauss' law, however, only a modified energy. Moreover, the equations become only explicit since the directions are decoupled by splitting the kinetic energy into the three subsystems of different directions. This separation is, however, limited to orthogonal grids and alternative time-stepping schemes are necessary once curvilinear coordinates are introduced. In this paper, we therefore study an alternative approach for the time discretization based on the form \eqref{eq:eqns_of_motion_Poisson} and the discrete gradient method.

\section{Implicit time stepping based on the discrete gradient method}\label{sec:disgrad_general}

\subsection{Discrete gradient time stepping}

The discrete gradient method was proposed by McLachlan et al.~\cite{McLachlan99} as a general method to construct energy conserving time stepping for conservative PDEs in skew-symmetric form, i.e. for a semi-discretization of the form
\begin{equation*}
 \dot{\ub} = \mathcal{J}(\ub) \cdot \DD_{\ub} \hat{\mathcal{H}} ( \ub) \quad \text{with} \quad \mathcal{J}(\ub)^\top = -\mathcal{J}(\ub).
\end{equation*}
The discrete gradient $\bar{\nabla} \mathcal{H}\left(\ub^{m}, \ub^{m+1}\right)$ for time step $[t_m,t_{m+1}]$ shall then satisfy
\begin{equation*}
\left( \ub^{m+1}-  \ub^m \right)^\top \bar{\nabla}  \hat{\ham}(\ub^m, \ub^{m+1}) =  \hat{\ham}(\ub^{m+1}) -  \hat{\ham}(\ub^m).
\end{equation*}
For any skew-symmetric approximation $\bar{\mathcal{J}}$ of $\mathcal{J}$, the following implicit scheme is then energy conserving,
$$ \frac{\ub^{m+1}-\ub^m}{\Delta t}  = \bar{\mathcal{J}} \bar{\nabla}  \hat{\ham}(\ub^m, \ub^{m+1}).$$
Energy conservation can be easily seen by the following calculation,
\begin{align*}		
		\hat{\ham}(\ub^{m+1}) -  \hat{\ham}(\ub^m) &= \left( \ub^{m+1} - \ub^m \right)^\top \bar{\nabla}  \hat{\ham}(u^m, u^{m+1}) \\&= \Delta t \bar{\nabla}  \hat{\ham}(\ub^m, \ub^{m+1})^{\top} \bar{\mathcal J}^{\top} \bar{\nabla}  \hat{\ham}(\ub^m, \ub^{m+1}) \\&= -\Delta t \bar{\nabla}  \hat{\ham}(\ub^m, \ub^{m+1})^{\top} \bar{\mathcal J} \bar{\nabla}  \hat{\ham}(\ub^m, \ub^{m+1}) = 0.
		\end{align*}
There is some freedom in the choice of the discrete gradient method. One systematic way of constructing a discrete gradient method is the average vector field method \cite{Celledoni12} that defines the time step as 
\begin{equation}\label{eq:avf}
 \frac{\ub_{m+1}-\ub_m}{\Delta t} = \int_0^1 g((1-\xi)\ub_m+\xi \ub_{m+1}) \du \xi, 
\end{equation}
where $g(\ub) = \mathcal{J}(\ub)D_{\ub} \hat{\ham}(\ub)$.

\subsection{Discrete gradients and antisymmetric splitting of the Poisson matrix}

The discrete gradient method applied to the full Vlasov--Maxwell system yields a heavily non-linear system for the time discretization. In order to simplify the system, we may split the Poisson matrix into several antisymmetric submatrices. Then, we can apply the discrete gradient method separately to each subsystem and combine the solutions of the subsystems in a splitting method, e.g. a Lie splitting for first order, Strang splitting for second order, or composition methods for higher order (see \cite{HairerLubichWanner:2006} for a splitting methods of various orders). 

The Poisson matrix \eqref{eq:poisson_matrix} can for instance be split into the following four antisymmetric matrices:
{
\begin{subequations}\label{eq:poisson_matrix_splitting}
\begin{align}
J_1 &:=\begin{pmatrix}
0 & {\MM_p^{-1}} & 0 & 0 \\
{- \MM_p^{-1}} & 0 & 0 \\
0 & 0 & 0 &  0\\
0 & 0 & 0 & 0 \\
\end{pmatrix} .\\
J_2 &:=\begin{pmatrix}
0 & 0 & 0 & 0 \\
0 & {\MM_p^{-1} \MM_q \BB (\X,\cb) \, \MM_p^{-1}} & 0 & 0 \\
0 & 0 & 0 & 0 \\
0 & 0 & 0 & 0 \\
\end{pmatrix} .\\
J_3 &:=\begin{pmatrix}
0 & 0 & 0 & 0 \\
0 &0 & 0 & 0 \\
0 & 0 & 0 &  {\MM_{1}^{-1} \C^\top}\\
0 & 0 & {- \C \MM_{1}^{-1}} & 0 \\
\end{pmatrix} .\\
J_4 &:=\begin{pmatrix}
0 & 0 & 0 & 0 \\
0 & 0 & {\MM_p^{-1} \MM_q \LaB^1(\X) \MM_{1}^{-1}} & 0 \\
0 & {- \MM_{1}^{-1} \LaB^1 (\X)^\top \MM_q \MM_p^{-1}} & 0 &  0\\
0 & 0 & 0 & 0 \\
\end{pmatrix} .
\end{align}
\end{subequations}
}
This yields the following four subsystems,
\begin{enumerate}
\item System 1: { $\dot{\X}= \V$}.
\item System 2: {$\dot{\V}= \MM_p^{-1} \MM_q \BB (\X,\cb) \V$}.
\item System 3: {$\dot{\ce}= \MM_{1}^{-1} \C^\top \MM_{2} \cb$, $\dot{\cb}= - \C \ce$}.
\item System 4: {$\dot{\V}= \MM_p^{-1} \MM_q \LaB^1 (\X) \ce $, $\dot{\ce}= -M_{1}^{-1}\LaB^1 (\X)^\top \MM_q \V$}.
\end{enumerate}
The first two systems can be solved analytically while we need to define a suitable discrete gradient for the last two systems. Below, we give the analytic solutions of systems 1 and 2. Discrete gradients for system 3 and 4 (possibly combined with 1) are the subject of Section \ref{sec:avf} and Section \ref{sec:disgrad}.

\subsubsection{Solution of system 1}

In the first subsystem, only the position is updated and the right-hand-side is independent of $\X$. The solution advancing the equation from time $t_0$ to $t$ is hence given as
\begin{equation*}
\X(t) = \X(t_0) + (t-t_0) \V.
\end{equation*} 

\subsubsection{Solution of system 2}

The second system only updates the velocity and it decomposes into one equation for each particle, namely
\begin{equation*}
\frac{\du}{\du t}\vb_a = \frac{q_s}{m_s} \widehat{\Bb}_{h} (\xb_a) \vb_a =  \frac{q_s}{m_s}\vb_a\times \Bb(\xb_a) , \quad a=1, \ldots, N_p.
\end{equation*}
This is a rotation round the magnetic axis $\tilde{\bb} = \frac{1}{ \|\Bb(\xb_a)\|_2}(B_{h,1}(\xb_a),B_{h,2}(\xb_a) , B_{h,3}(\xb_a))^\top$, denoting
by $ \|\Bb(\xb_a)\|_2=\sqrt{B_{h,1}(\xb_a)^2+B_{h,2}(\xb_a)^2 + B_{h,3}(\xb_a)^2}$, (that is fixed over time in this subsystem) with an angle  $\alpha = \Delta t \frac{q_s}{m_s} |\Bb(\xb_a)|$ depending on the time step $\Delta t = t-t_0$. The solution is given by
\begin{equation*}
\vb_a(t) = R(\Bb_h(\xb_a)) \vb_a(t_0), 
\end{equation*}
where $R(\Bb_h(\xb_a))$ is the rotation matrix 
{\footnotesize
\begin{eqnarray}\label{eq:rotation_matrix}
R(\Bb_h(\xb_a))=\begin{pmatrix}
\tilde{b_1}^2 + \left(\tilde{b}_2^2+\tilde{b}_3^2\right) \cos(\alpha) & \tilde{b}_3 \sin(\alpha)+ \tilde{b}_2\tilde{b}_1 (1-\cos(\alpha)) &-\tilde{b}_2 \sin(\alpha)+ \tilde{b}_3\tilde{b}_1 (1-\cos(\alpha)) \\
-\tilde{b}_3 \sin(\alpha)+ \tilde{b}_2\tilde{b}_1 (1-\cos(\alpha)) & \tilde{b_2}^2 + \left(\tilde{b}_1^2+\tilde{b}_3^2\right) \cos(\alpha) &\tilde{b}_1 \sin(\alpha)+ \tilde{b}_3\tilde{b}_2 (1-\cos(\alpha))  \\
\tilde{b}_2 \sin(\alpha)+ \tilde{b}_3\tilde{b}_1 (1-\cos(\alpha)) &-\tilde{b}_1 \sin(\alpha)+ \tilde{b}_3\tilde{b}_2 (1-\cos(\alpha)) & \tilde{b_3}^2 + \left(\tilde{b}_2^2+\tilde{b}_1^2\right) \cos(\alpha)
\end{pmatrix}.
\end{eqnarray}
}

\section{ A semi-discrete average vector field discretization}\label{sec:avf}

The systems 3 and 4 are both linearly implicit since the right-hand-sides are linearly dependent on a dynamic variable changing in the respective step. We now construct a time stepping for both systems that is based on the average vector field method \eqref{eq:avf}.

\subsection{Solution of System 3}\label{sec:avf_advect_eb}

Applying the average-vector-field method to the equations of System 2, we get the following system of linear equations for the unknown coefficients $\ce^{m+1}, \cb^{m+1}$ ,

\begin{equation*}
\begin{pmatrix}
\MM_{1} & -\frac{\Delta t}{2} \C^T \MM_{2} \\
\frac{\Delta t}{2} \MM_{2} \C &  \MM_{2}
\end{pmatrix} \begin{pmatrix}
\ce^{m+1} \\ \cb^{m+1}
\end{pmatrix} =\begin{pmatrix} \MM_{1} & \frac{\Delta t}{2} \C^T \MM_{2} \\
-\frac{\Delta t}{2} \MM_{2} \C &  \MM_{2}
\end{pmatrix} \begin{pmatrix}
\ce^{m} \\ \cb^{m}
\end{pmatrix}.
\end{equation*}
The system can, of course, be solved in this form. For this, an iterative GMRES solver can for instance be used. For increasing degree of the spline basis, the system, however, gets rather ill-conditioned and we need a good preconditioner to solve the system. A simple preconditioner would be to split the equations into two explicit equations, i.e. to first solve for $\cb^{m+1}$ for given $\ce^{m}$ and then to solve for $\ce^{m+1}$ (or vice versa).

On the other hand, the equations for $\ce^{m+1}$ and $\cb^{m+1}$ can be decoupled using the Schur complement $S=\MM_{1} + \frac{\Delta t^2}{4} \C^T \MM_{2} \C$:
\begin{equation*}
\begin{pmatrix}
\MM_{1} & -\frac{\Delta t}{2} \C^T \MM_{2} \\
\frac{\Delta t}{2} \MM_{2} \C &  \MM_{2}
\end{pmatrix}^{-1} = \begin{pmatrix}
I & 0 \\ -\frac{\Delta t}{2} \C & I
\end{pmatrix} \begin{pmatrix}
S^{-1} & 0 \\ 0 & \MM_{2}^{-1}
\end{pmatrix} \begin{pmatrix}
I & \frac{\Delta t}{2} \C^T\\ 0 & I
\end{pmatrix}.
\end{equation*}
With this expression for the matrix inverse, we get the following two equations:
\begin{subequations}
\begin{align}
\ce^{m+1} &= S^{-1} \left( \left( \MM_{1}- \frac{\Delta t^2}{4} \C^T \MM_{2} \C\right) \ce^m + \Delta t \C^T \MM_{2} \cb^m \right), \label{eq:system3_e} \\
\cb^{m+1} &= \cb^{m} - \frac{\Delta t}{2} \C \left(\ce^m+\ce^{m+1}\right). \label{eq:system3_b}
\end{align}
\end{subequations}
Hence, we only need to solve the system $S \ce^{m+1} = \cf$ for $\ce^{m+1}$ with given right-hand-side $\cf$ and the magnetic field can then be updated by an explicit equation.

Let us now consider the implicit equation a bit more in detail. For this, we split the equation into three parts for each of the components of the field. The discrete mass matrix has block-diagonal form
\begin{equation*}
\MM_{2} = \begin{pmatrix}
M_{21} & 0 & 0 \\
0 & M_{22} & 0 \\
0 & 0 & M_{23}
\end{pmatrix}
\end{equation*}
and the discrete curl matrix has the following block structure
\begin{equation*}
\begin{pmatrix}
0 & -D_3 & D_2 \\
D_3 & 0 & -D_1 \\
-D_2 & D_1 & 0
\end{pmatrix},
\end{equation*}
where $D_i$, $i=1,2,3$, denotes the derivative matrix along direction $i$. With this notation, we have the following expression for $\C^\top \MM_{2} \C$:
\begin{eqnarray*}
\C^\top \MM_{2} \C = \begin{pmatrix}
D_3^\top M_{22} D_3 + D_2^\top M_{23} D_2 & -D_2^\top M_{23} D_1 & - D_3^\top M_{22} D_1 \\
-D_1^\top M_{23} D_2 &  D_3^\top M_{21} D_3 + D_1^\top M_{23} D_1 & -D_3^\top M_{21} D_2 \\
-D_1^\top M_{22} D_3 &- D_2^\top M_{21} D_3 & D_2^\top M_{21} D_2 + D_1^\top M_{22} D_1 \\
\end{pmatrix}.
\end{eqnarray*}
Componentwise the equation therefore reads
\begin{eqnarray*}
&\left( M_{11} + \frac{\Delta t^2}{2} \left( D_3^T M_{22} D_3 + D_2^T M_{23} D_2 \right) \right) \ce_1 - \frac{\Delta t^2}{2} D_2^T M_{23} D_1 \ce_2 - \frac{\Delta t^2}{2}D_3^T M_{22} D_1 \ce_3 = \cf_1,\\
&-\frac{\Delta t^2}{2} D_1^T M_{23} D_2 \ce_1 + \left( M_{12} + \frac{\Delta t^2}{2}\left( D_3^T M_{21} D_3 + D_1^T M_{23} D_1 \right) \right) \ce_2 -\frac{\Delta t^2}{2}D_3^T M_{21} D_2 \ce_3 = \cf_2, \\ 
&-\frac{\Delta t^2}{2}D_1^T M_{22} D_3 \ce_1 - \frac{\Delta t^2}{2}D_2^T M_{21} D_3 \ce_2 +\left( M_{13} +\frac{\Delta t^2}{2} \left( D_2^T M_{21} D_2 + D_1^T M_{22} D_1\right)\right)\ce_3 = \cf_3.
\end{eqnarray*}

\subsubsection{Direct inversion in Fourier space for periodic boundary conditions}

In this paper, we consider the special case of a periodic box. Then, all the one-dimensional matrices are circulant and can thus be  diagnonalized by Fourier transformation. The derivative matrix is given by the circulant matrix 
\begin{equation*}
 D = \frac{1}{\Delta x}\begin{pmatrix}
1 & && -1\\
-1 & 1 & &\\
& \ddots & \ddots & \\
 & & -1 & 1
\end{pmatrix}.
\end{equation*} 
The eigenvalues of the one dimensional building blocks of the matrix $\C^\top \MM_{2} \C$ can therefore computed to be
\begin{itemize}
	\item $D$: $\lambda_k^{+} = \frac{1}{\Delta x} \left( 1- \exp \left(-\frac{2\pi \im k}{n} \right) \right)$, $j=0, \ldots, n-1$. 
	\item $D^T$: $\lambda_k^{-} = \frac{1}{\Delta x} \left( 1- \exp \left(\frac{2\pi \im k}{n} \right) \right)$, $j=0, \ldots, n-1$.
	\item $M$ with row $c_p, \ldots, c_0, \ldots, c_p$ (p order of the spline): $\lambda_k^{(p)} = c_0 +\sum_{j=1}^p c_j 2 \cos\left( \frac{2\pi k j}{n}\right)$.
\end{itemize}
After Fourier transformation, we end up with a $3\times 3$ system for each Fourier mode which can be solved explicitly.


\subsection{Solution of System 4}\label{sec:avf_advect_e}

The equations of System 4 can be split into three separate equations for the pairs $(\V_i,\ce_i)$, $i=1,2,3$, 
\begin{equation}\label{eq:system3_ode}
\dot{\V}_i= \MM_p^{-1} \MM_q \LaB^1_i (\X) \ce_i , \quad \dot{\ce}_i= -\MM_{1}^{-1}\LaB^1_i (\X)^\top \MM_q \V_i.
\end{equation}
Applying the average-vector-field method to equation \eqref{eq:system3_ode}, we get the following linear system
\begin{multline*}
\begin{pmatrix}
I & - \frac{\Delta t}{2} M_p^{-1} M_q \LaB^1_i(\X)  \\ \frac{\Delta t}{2} \LaB^1_i(\X)^\top M_q  & \MM_{1}
\end{pmatrix} \begin{pmatrix}
\V_i^{m+1} \\ \ce_i^{m+1}
\end{pmatrix} =\\
 \begin{pmatrix}
I &  \frac{\Delta t}{2} M_p^{-1} M_q \LaB^1_i(\X)  \\ -\frac{\Delta t}{2} \LaB^1_i(\X)^\top M_q  & \MM_{1}
\end{pmatrix} \begin{pmatrix}
\V_i^{m} \\ \ce_i^{m}
\end{pmatrix}.
\end{multline*}
Defining the Schur complement $S = \MM_{1} + \frac{\Delta t^2}{4} \MM_p^{-1}\MM_q\LaB^1_i(\X)^T M_q \LaB^1_i(\X)$, we get the following expression for the inverse of the left-hand-side matrix
\begin{multline*}
\begin{pmatrix}
I & - \frac{\Delta t}{2} M_p^{-1} M_q \LaB^1_i(\X)  \\ \frac{\Delta t}{2} \LaB^1_i(\X)^\top M_q  & \MM_{1}
\end{pmatrix} ^{-1} \\
= \begin{pmatrix}
I &  \frac{\Delta t}{2} M_p^{-1} M_q \LaB^1_i(\X)  \\ 0 & I
\end{pmatrix}  \begin{pmatrix}
I & 0  \\ 0 & S^{-1}
\end{pmatrix} \begin{pmatrix}
I & 0  \\ -\frac{\Delta t}{2} \LaB^1_i(\X)^\top M_q  & I
\end{pmatrix}.
\end{multline*}
Hence, the system can be solved in three steps:
\begin{enumerate}
	\item $v_i^* = v_i^m + \frac{\Delta t}{2} \MM_p^{-1}\MM_q \LaB^1_i(\X) \ce_i^m$,\\$\ce_i^* = \left( \MM_{1} - \frac{\Delta t^2}{4} \MM_p^{-1}\MM_q \LaB^1_i(\X)^T M_p \LaB^1_i(\X) \right) \ce_i^m - \Delta t \LaB^1( \X)^T M_q v_i^m$.
	\item $\ce_i^{m+1} = S^{-1} \ce_i^*$.
	\item $v_i^{m+1} = v_i^* + \frac{\Delta t}{2} \MM_p^{-1}\MM_q \LaB^1_i(\X) \ce_i^{m+1}$.
\end{enumerate}
Note that the implicit part (step 2) is reduced to the field equations and the particle equations can be solved explicitly.

\subsubsection{Preconditioning the linear solver}

Since the matrix $S$ is a symmetric matrix, we can solve the system in step 2 using the conjugate gradient method. However, also $S$ is increasingly ill-conditioned for higher degree of the splines so that we need a good preconditioner. 
Let us consider the matrix $N_i:=\LaB^1_i(\X)^T M_q \LaB^1_i(\X)$ with elements
\begin{eqnarray*}
(N_i)_{j,k} = q \sum_a w_a \Lambda_{i,j}(\xb_a) \Lambda_{i,k}(\xb_a). 
\end{eqnarray*}
We note that this is a Monte Carlo approximation of the integral $\int_{\Omega} \rho(\xb) \Lambda_{i,j}(\xb) \Lambda_{i,k}(\xb) \du \xb$. If the charge density was equal to one over the spatial domain, this would be an approximation of the mass matrix. Therefore, we refer to this matrix as \textit{particle sampled mass matrix}. Since the perturbations from equilibrium are usually small, $\rho \ll 1$ holds true so that the particle sampled mass matrix is much smaller than one and $\MM_{1}$ is a reasonable approximation for $S$ and can be used as a preconditioner. This approximation is the better the smaller the time step, the more particles per cell are used, and the smaller the perturbations. For the periodic box, the mass matrix is circulant and can be inverted in Fourier space (cf. Sec.~\ref{sec:avf_advect_eb}). 

Note that for a simulation of electrons with a neutralizing electron background, $N$ only contains the electron background. In this case, the sampled function is $1-\rho$ and, hence, close to one for small perturbations. In this case, $N$ is close to a Monte-Carlo approximation of the mass matrix and $(1+\frac{\Delta t^2}{4}\frac{q^2}{m} )\MM_{1}$ is a good approximation to the matrix $S$.


\subsection{Summary of the average-vector-field time stepping}

The proposed semi-implicit and energy-conserving scheme, called average-vector-field scheme in the following, is composed of the four operators resulting from the splitting of the Poisson matrix \eqref{eq:poisson_matrix_splitting}. For the time-stepping, the following four operators are combined:
\begin{enumerate}
	\item Operator 1: $\X(t) = \X(t_0) + (t-t_0) \V$.
	\item Operator 2: $\vb_a(t) = R(\Bb_h(\xb_a)) \vb_a(t_0)$ with the rotation matrix defined by eqn.~\eqref{eq:rotation_matrix}.
	\item Operator 3:\begin{align}
\ce(t) &= S^{-1} \left( \left( \MM_{1}- \frac{\Delta t^2}{4} \C^T \MM_{2} \C\right) \ce(t_0) + \Delta t \C^T \MM_{2} \cb(t_0) \right), \\
\cb(t) &= \cb(t_0) - \frac{\Delta t}{2} \C \left(\ce(t_0)+\ce(t)\right). 
\end{align}
\item Operator 4: \begin{enumerate}
	\item $v_i^* = v_i(t_0) + \frac{\Delta t}{2} \MM_p^{-1}\MM_q \LaB^1_i(\X) \ce_i(t_0)$,\\$\ce_i^* = \left( \MM_{1} - \frac{\Delta t^2}{4} \MM_p^{-1}\MM_q \LaB^1_i(\X)^T M_p \LaB^1_i(\X) \right) \ce_i(t_0) - \Delta t \LaB^1( \X)^T M_q v_i(t_0)$.
	\item $\ce_i(t) = S^{-1} \ce_i^*$,  $S = \MM_{1} + \frac{\Delta t^2}{4} \MM_p^{-1}\MM_q\LaB^1_i(\X)^T M_q \LaB^1_i(\X)$.
	\item $v_i(t) = v_i^* + \frac{\Delta t}{2} \MM_p^{-1}\MM_q \LaB^1_i(\X) \ce_i(t)$.
\end{enumerate}
\end{enumerate}
The four operators can now be combined in various ways to build the full time step. The first order Lie splitting is build upon full time steps of each of the operators one after the other. Second order can be obtained when combining two Lie splitting steps with opposite ordering of the operators of half a time step each. In this case, the last operator of the first Lie step and the first operator of the second Lie step are the same so that we can instead place a full time step of one operator in the middle. Clearly, for a large number of particles per cell, the fourth operator is the most expensive one. Therefore, the complexity of the algorithm is reduced when placing this operator in the middle and call it only once. 

A Strang splitting of the following form yields a shortest run time: half time step of operator 3, half time step of operator 1, half time step of operator 2, full time step with operator 4, half time step with operator 2, half time step with operator 1, half time step with operator 3. With this ordering we only have to apply the most expensive operator 4 once and we can merge the updates of operator 1 and 2 (that do not touch the fields) into the particle loops of operator 4. The number of loops over all particles is minimized in order to maximize the arithmetic intensity of the algorithm. The computational complexity is then dominated by the assembling of the particle sampled mass matrix. This results in the algorithm outlined in Algorithm \ref{alg:vm_avf}. We note that the ordering of operator 1 and 2 could also be exchanged without changing the structure of the loops. 

On the other hand, we can also ask which ordering gives best accuracy. In our examples, we found that better accuracy is sometimes achieved when placing operator 3 after the operators 1 and 2. However, in this case we need to traverse all particles two times more per time step (or one time if we fuse this step between two time steps in a ``first same as last'' procedure). 

The computational overhead compared to the explicit Hamiltonian splitting is rather limited, namely we need to solve the curl-part of Maxwell's equation implicitly (which is, however, a small problem compared to the particle loops) and we need the particle sampled mass matrix. On the other hand, no integrals need to be computed for the deposition of the curl.


The method (with any ordering of the operators) is energy-conserving since the single operators are solved either exactly or based on an energy-conserving average-vector-field discretization. On the other hand, the Casimir invariants can be destroyed by the splitting of the Poisson matrix. In particular, Gauss' law is not conserved over time with this time stepping method as will be discussed in the next section.

\begin{algorithm}
\caption{Second order average-vector-field scheme.}\label{alg:vm_avf}
\For{$m=1, \ldots$}{
$\ce^{*} = S^{-1} \left( \left( \MM_{1}- \frac{\Delta t^2}{16} \C^T \MM_{2} \C\right) \ce^m + \Delta t \C^T \MM_{2} \cb^m \right)$ \;
$\cb^{*} = \cb^{m} - \frac{\Delta t}{4} \C \left(\ce^m+\ce^{*}\right)$\;
$\cb = 0$\;
$N = 0$\;
\For{$a=1, \ldots, N_p$}{
	$\xb_a^* = \xb_a^m + \frac{\Delta t}{2} \vb_a^m$\;
	$\vb_a^* = R(\bb^*) \vb_a^m$\;
	$\vb_a^{**} = \vb^* + \frac{\Delta t}{2} \frac{q_s}{m_s} \LaB^1(\xb_a^*) \eb^*$\;
	$\cj = \cj + \Delta t q_s w_a \LaB^1(\xb_a^*) \vb_a$\;
	\For{$i=1,2,3$}{
	$N_i = N_i + \frac{\Delta t^2}{4} q_s w_a  \LaB_i^1(\xb_a^*) \LaB_i^1(\xb_a^*)$ \;
	}
}
$\ce^{**} = \left(\MM_{1} - \begin{pmatrix}
N_1 & 0 & 0 \\ 0 & N_2 & 0 \\ 0 & 0 & N_3
\end{pmatrix}\right) \ce^* - \cj$\;
\For{$a=1, \ldots, N_p$}{
	$\vb_a^{***} = \vb^{**} + \frac{\Delta t}{2} \frac{q_s}{m_s} \LaB^1(\xb_a^*) \eb^{**}$\;
	$\vb_a^{m+1} = R(\bb^*) \vb_a^{***}$\;
	$\xb_a^{m+1} = \xb_a^{*} + \frac{\Delta t}{2} \vb_a^{m+1}$\;
}
$\ce^{m+1} = S^{-1} \left( \left( \MM_{1}- \frac{\Delta t^2}{16} \C^T \MM_{2} \C\right) \ce^{**} + \Delta t \C^T \MM_{2} \cb^* \right)$ \;
$\cb^{m+1} = \cb^{*} - \frac{\Delta t}{4} \C \left(\ce^{**}+\ce^{m+1}\right)$\;
}
\end{algorithm}


\section{Conservative implicit discrete gradient method}\label{sec:disgrad}

In the previous section, we have derived an energy-conserving semi-implicit time propagator. However, it was shown in \cite[Sec. 4.6]{Kraus17} that the semi-discrete equations of motion of the GEMPIC framework (cf. \eqref{eq:vlasov_maxwell_equations_of_motion}) satisfy Gauss' law which the average-vector-field scheme does not as we will show in the following. In this section, we derive an alternative discrete gradient method that---in addition to energy---also conserves Gauss' law.


\subsection{Average-vector-field method destroys Gauss' law}

In order to satisfy Gauss' law, it is important how the current is accumulated in Amp\`ere's law,
\begin{equation}\label{eq:gauss_essential}
\MM_{1} \ce^{m+1} = \MM_{1} \ce^m - \int_{t_m}^{t_{m+1}} \LaB^1 (\X(\tau))^\top \MM_q \V(\tau) \du \tau.
\end{equation}
Applying $\G^\top$ to \eqref{eq:gauss_essential} yields,
\begin{equation*}\begin{aligned}
\G^\top \MM_{1} \ce^{m+1} &= \G^\top \MM_{1} \ce^m - \int_{t_m}^{t_{m+1}} \G^\top\LaB^1 (\X(\tau))^\top \MM_q \V(\tau) \du \tau .
\end{aligned}\end{equation*}
If $\frac{\du \X(\tau)}{\du \tau} = \V(\tau)$ holds true, we can use the chain rule to identify the integrand as a time derivative, namely
\begin{equation*}
\G^\top\LaB^1 (\X(\tau))^\top \MM_q \V(\tau) = \frac{\text{d}\phantom{\tau} }{\text{d} \tau}  \LaB^1(\X(\tau))^\top \MM_q \mathbb{1}_{N_p}.
\end{equation*}
Hence, the integral over $t$ can be evaluated as
\begin{equation*}\begin{aligned}
\G^\top \MM_{1} \ce^{m+1} &= \G^\top \MM_{1} \ce^m - \int_{t_m}^{t_{m+1}} \G^\top\LaB^1 (\X(\tau))^\top \MM_q \V(\tau) \du \tau \\ 
&= \G^\top \MM_{1} \ce^m -\LaB^1(\X(t_{m+1}))^\top \MM_q \mathbb{1}_{N_p} + \LaB^1(\X(t_{m}))^\top \MM_q \mathbb{1}_{N_p}.
\end{aligned}\end{equation*}
This means that, if the discrete version of Gauss' law $\G^\top \MM_{1} \ce^m = \LaB^1(\X(t_{m}))^\top \MM_q \mathbb{1}_{N_p}$ holds at time $t_m$, it also holds at time $t_{m+1}$.  

From this derivation, it becomes clear that the $\X$ update and the current deposition should not be split. Hence, System 1 and System 4 need to be merged to be able to conserve Gauss' law. However, the average-vector-field ansatz itself already destroys Gauss' law: An average-vector-field discretization of the combined System 1 and 4 reads,
\begin{eqnarray*}
\X^{m+1} &=& \X^m + \frac{\Delta t}{2} \left( \V^{m+1} + \V^m \right), \\
\MM_{1} \ce^{m+1} &=& \MM_{1} \ce^m - \\
&&\Delta t\int_0^1 \LaB^1 \left(\X^m + \xi \frac{\Delta t}{2 }\left( \V^{m+1} + \V^m \right)\right)^\top \MM_q \left(\V^{m+1} \xi + \V^m (1-\xi)\right) \du \xi.
\end{eqnarray*}
In this case, Gauss's law would be conserved if $\frac{\du  \left(\X^m + \xi \frac{\Delta t}{2 }\left( \V^{m+1} + \V^m \right)\right)}{\du \xi} \stackrel{!}{=}    \Delta t\V^{m+1} \xi + (1-\xi) \V^m$. However, it holds that $\frac{\du  \left(\X^m + \xi \frac{\Delta t}{2 }\left( \V^{m+1} + \V^m \right)\right)}{\du \xi} =   \frac{\Delta t}{2} \left(\V^{m+1} + \V^m \right)$ in this scheme. 

We note that the curl-part of Amp\`ere's law does not change the divergence of the electric field since the discrete gradient and curl operators respect the relation $\C  \G = 0$: Multiplying the curl-part of Amp\`ere's law by $\G^\top$, we get
\begin{equation*}
\G^\top \MM_{1} \dot \eb - \frac{\Delta t}{2} \G^\top \C^\top  \MM_{2} \bb = 0,
\end{equation*}
with the second term being zero since $\G^\top \C^\top = (\C \G)^\top = 0$. Hence the divergence of the electric field does not change with time.
In particular, the divergence is not changed in the solution of System 3, where the electric field is updated by equation \eqref{eq:system3_e}. Multiplying the equation by $\G^\top$ yields
\begin{equation*}
\G^\top \left( \MM_{1} + \frac{\Delta t^2}{4}\C^\top \MM_{2} \C \right) \eb^{m+1} = \G^\top \left( \MM_{1} - \frac{\Delta t^2}{4} \C^\top \MM_{2} \C^\top \right) \eb^m + \Delta t \G^\top \C^\top \MM_{2} \bb^{m}.
\end{equation*}
Again exploting $\G^\top \C^\top = (\C \G)^\top = 0$ yields
\begin{equation*}
\G^\top \MM_{1} \eb^{m+1} = \G^\top \MM_{1}  \eb^m.
\end{equation*}

\subsection{Alternative discrete gradient that conserves Gauss' law}

Nevertheless, a discrete gradient method that conserves Gauss' law can be constructed by not only combining system 1 and 4 but also modifying the definition of the discrete gradient. For the new construction of the discrete gradient, we assume that
the particle trajectory is linear between $\X^m$ and $\X^{m+1}$ with velocity $(\V^m+\V^{m+1})/2$ and that the electric field is constant with its average value, so that on this interval $\V(\tau)= (\V^m+\V^{m+1})/2$ is constant and the approximate trajectory is defined by 
$$\X(\tau) =  ((t_{m+1}-\tau)\X^m + (\tau-t_m) \X^{m+1})/\Delta t $$ and 
$$\frac{\du \X}{\du \tau} = \frac{\X^{m+1}-\X^m}{\Delta t}=\frac{\V^m+\V^{m+1}}{2}.$$
The discrete gradient method with this approximation leads to
\begin{subequations}\label{eq:disgrad}\begin{align}
\frac{\X^{m+1}-\X^m}{\Delta t}&=\frac{\V^m+\V^{m+1}}{2}, \\
\frac{\V^{m+1}-\V^m}{\Delta t}&= \MM_p^{-1} \MM_q \frac {1}{\Delta t} \int_{t_m}^{t_{m+1}}\LaB^1 (\X(\tau))\du\tau \left( \frac{\ce^m+\ce^{m+1}}{2} \right),\label{eq:disgrad_v}\\
\frac{\ce^{m+1}-\ce^m}{\Delta t} &= -\MM_{1}^{-1}\frac {1}{\Delta t}  \int_{t_m}^{t_{m+1}} \LaB^1 (\X(\tau))^\top \du\tau  \MM_q \left(\frac{\V^m+\V^{m+1}}{2}  \right),\label{eq:disgrad_e}
\end{align}\end{subequations}
with $\X(\tau) = \X^m\frac{t_{m+1}-t}{t_{m+1}-t_m} + \X^{m+1}\frac{t-t_m}{t_{m+1}-t_m}$. The nonlinear system \eqref{eq:disgrad} can be reformulated as a fixed point iteration in $\U= \begin{pmatrix}
\X^{m+1} \\
\V^{m+1} \\
\eb^{m+1}
\end{pmatrix}$.

It can easily be shown that the scheme \eqref{eq:disgrad} is indeed energy conserving by the following calculations
\begin{eqnarray*}
&&\left(\V^{m+1}\right)^\top\MM_p\V^{m+1}-\left(\V^m\right)^\top \MM_p \V^m\\
&&\stackrel{\eqref{eq:disgrad_v}}{=}  \left( (\V^{m+1} + \V^m)^\top
\MM_q  \int_{t_m}^{t_{m+1}}\LaB^1_i (\X(\tau))\du\tau \right) \left( \frac{\ce^m+\ce^{m+1}}{2} \right) \\
&&\stackrel{\eqref{eq:disgrad_e}}{=} -\left(\left(\ce^{m+1}\right)^\top\MM_1\ce^{m+1}-\left(\ce^m\right)^\top \MM_1 \ce^m\right).
\end{eqnarray*}
Hence, the difference in the kinetic energy equals the negative difference in the potential energy (since the magnetic field is not changed in this system) and thus the total energy is conserved.
\subsubsection{Summary of the discrete gradient time stepping}
Scheme \eqref{eq:disgrad} can be combined with the analytic solution of System 2 and the average-vector-field solution of System 3 (as discussed in Section \ref{sec:avf_advect_eb}) to an implicit propagator that conserves both energy and Gauss' law. 

System \eqref{eq:disgrad} is now not only linearly implicit (as it was the case when the $\X$ propagation was separated from $\V$ and $\ce$) so that we need to use an iterative algorithm to solve this nonlinear equation. In our implementation, we use the average-vector-field scheme to produce a good starting point for the nonlinear iteration and then continue to improve the approximation with Picard iteration until a predefined tolerance is met. Note that this scheme also comes with a memory overhead since we need to store a second particle position and velocity for all particles during the nonlinear iteration. Algorithm \ref{alg:vm_dg} summarizes the proposed conservative discrete gradient scheme in a second order Strang splitting with minimal number of particle loops. However, we note that we need one particle loop per Picard iteration in this case. Reordering operator 2 and 3 and spending one (or two) extra particle loops per time step will only marginally increase the computing time for each time step and can therefore easily pay-off if accuracy is increased.

\subsubsection{Exact numerical evaluation of the integrals}
Note that we have to evaluate the integrals of the form $\int_0^1\Lambda_j^1(\X(\tau)) \du \tau $, where $\Lambda_j^1$ is a spline of a certain degree in each of the three variables. As a function of $t$ the integrand is then a polynomial of degree $p_1p_2p_3$ (where $p_i$ denotes the degree of the spline in each direction $i=1,2,3$) locally in each cell of the three dimensional domain. In order to solve the integral exacly (which is necessary to conserve Gauss' law), we can use Gauss--Legendre quadrature with $\lceil \frac{p_1p_2p_3+1}{2}\rceil$ points separately in each cell crossed by the line integral. Note that the degree of the quadrature is generally quite high, on the other hand, we only need a one-dimensional quadrature rule. As an example, if we use splines of degree 3 for the 0-form, i.e. we have two times degree 3 and once degree 2 in the 1-form, we need to use a 10-point-quadrature. For splines of degree 2, on the other hand, we only need 3 points, so the complexity rapidely increases with the degree of the splines.

\begin{algorithm}
\caption{Second order conservative discrete gradient scheme.}\label{alg:vm_dg}
\For{$m=1, \ldots$}{
$\ce^{*} = S^{-1} \left( \left( \MM_{1}- \frac{\Delta t^2}{16} \C^T \MM_{2} \C\right) \ce^m + \Delta t \C^T \MM_{2} \cb^m \right)$ \;
$\cb^{*} = \cb^{m} - \frac{\Delta t}{4} \C \left(\ce^m+\ce^{*}\right)$\;
$\cj = 0$\;
$N = 0$\;
\For{$a=1, \ldots, N_p$}{
	$\vb_a^* = R(\bb^*) \vb_a^m$\;	
	$\xb_a^{\ii} = \xb_a^m + \frac{\Delta t}{2} \vb_a^*$\;
	$\vb_a^{**} = \vb^* + \frac{\Delta t}{2} \frac{q_s}{m_s} \LaB^1(\xb_a^{\ii}) \eb^*$\;
	$\cj = \cj + \Delta t q_s w_a \LaB^1(\xb_a^{\ii}) \vb_a$\;
	\For{$i=1,2,3$}{
	$N_i = N_i + \frac{\Delta t^2}{4} q_s w_a  \LaB_i^1(\xb_a^{\ii}) \LaB_i^1(\xb_a^{\ii})$ \;
	}
}
$\ce^{\ii} = \left(\MM_{1} - \begin{pmatrix}
N_1 & 0 & 0 \\ 0 & N_2 & 0 \\ 0 & 0 & N_3
\end{pmatrix}\right) \ce^* - \cj$\;
\For{$a=1, \ldots, N_p$}{
	$\vb_a^{\ii} = \vb_a^{**} + \frac{\Delta t}{2} \frac{q_s}{m_s} \LaB^1(\xb_a^{\ii}) \eb^{\ii}$\;
	$\xb_a^{\ii} = \xb_a^{\ii} + \frac{\Delta t}{2} \vb_a^{\ii}$\;
}
$residual=tolerance+1$\;
\While{$residual<tolerance$}{
$\cj = 0$\;
\For{$a=1, \ldots, N_p$}{
	$\bar{\vb}_a =\frac{1}{2}\left( \vb_a^* + \vb^{\ii}\right)$\;
	$\xb_a^{\ii} = \xb_a^{*} + \Delta t \bar{\vb}_a$\;
	$\cj =  \cj + \Delta t q_s w_a \int_0^1\LaB^1(\xb_a^* + \tau \bar{\vb}_a) \du \tau \vb_a$\;
	$\vb_a^{\ii} = \vb^{*} + \frac{\Delta t}{2} \int_{0}^1\frac{q_s}{m_s} \LaB^1(\xb_a^* + \tau \bar{\vb}_a) \du \tau \eb^{**}$\;
}
$\ce^{\text{old}} = \ce^{\ii}$\;
$\ce^{\ii} = \ce^* - M^{-1}\cj$\;
$residual = \|\ce^{\ii}-\ce_{\text{old}}\|_2$\;
}
\For{$a=1, \ldots, N_p$}{
	$\vb_a^{m+1} = R(\bb^*) \vb_a^{\ii}$\;
	$\xb_a^{m+1} = \xb_a^{\ii} + \frac{\Delta t}{2} \vb_a^{m+1}$\;
}
$\ce^{m+1} = S^{-1} \left( \left( \MM_{1}- \frac{\Delta t^2}{16} \C^T \MM_{2} \C\right) \ce^{\ii} + \Delta t \C^T \MM_{2} \cb^* \right)$ \;
$\cb^{m+1} = \cb^{*} - \frac{\Delta t}{4} \C \left(\ce^{\ii}+\ce^{m+1}\right)$\;
}
\end{algorithm}

\subsection{Substepping}

Another useful feature for the accurate simulation of low frequency phenomena, where the fields are slowly varying in time but on the grid scale in space, would be particle subcycling. We therefore adapt the method proposed by Chen, Chacon, \& Barnes \cite{Chen11} to our context. For this, we split the time step into $N_{sub}$ subintervals $[t_{\nu},\tau_{\nu+1}]$, $\nu=0,\ldots, N_{sub}$, of not necessarily identical length $\Delta \tau_\nu=\tau_{\nu+1}-\tau_\nu$ with $\tau_0=t_m$ and $\tau_{N_\nu}=t_{m+1}$:  $\Delta t= t_{m+1} -t_m =\sum_{\nu=0,N_{sub}-1} (\tau_{\nu+1}-\tau_\nu) $.
Then, keeping the electric field constant over all the substeps, we push the particles according to
$\X_m^0=\X^m$, $\V_m^0=\V^m$ and for $\nu=0$ to $N_\nu-1$ 
\begin{subequations}\label{eq:disgrad_sub}\begin{align}
\frac{\X^{\nu+1}_m-\X^\nu_m}{\Delta \tau_\nu}&=\frac{\V^\nu_m+\V^{\nu+1}_m}{2} \\
\frac{\V^{\nu+1}_m-\V^\nu_m}{\Delta \tau_\nu}&= \MM_p^{-1} \MM_q \frac {1}{\Delta \tau_\nu} \int_{\tau_\nu}^{\tau_{\nu+1}}\LaB^1 (\X_m(\tau))\du\tau \left( \frac{\ce^m+\ce^{m+1}}{2} \right)\label{eq:disgrad_sub_v}
\end{align} \end{subequations}
with $\X_m(\tau) =  ((\tau_{\nu+1}-\tau) \X^\nu_m + (\tau-\tau_\nu)  \X^{\nu+1}_m)/\Delta \tau_\nu$. The current for the update of the electric field then also needs to be updated as a sum over the contributions of each substep following the discrete gradient formulation in \eqref{eq:disgrad_e} by
\begin{equation}\label{eq:dgsub_e}
\frac{\ce^{m+1}-\ce^m}{\Delta t} =  -M_{1}^{-1}\frac {1}{\Delta t} 
\sum_{\nu=0}^{N_{sub}-1} \left(  \int_{\tau_\nu}^{\tau_{\nu+1}}\LaB^1 (\X_m(\tau))^\top\du\tau \right) \MM_q \left( \frac{\V^\nu_m+\V^{\nu+1}_m}{2} \right).
\end{equation}
For this version with substepping, we have two nested nonlinear systems. The sub-iteration \eqref{eq:disgrad_sub} is a nonlinear system in $\left(X^{\nu+1}_m,V^ {\nu+1}_m\right)$ which decomposes into separate systems for each particle which can be included into the overall iteration that is similar to the algorithm without substepping. In our implementation, we use Picard iterations for both nonlinear systems. A Strang splitting version with substepping is shown in Algorithm \ref{alg:vm_dgsub}.

\subsubsection{Conservation properties}

Also with substepping, the scheme is energy conserving as can be seen from the following calculation: Setting $\V^{m+1}=\V^{N_\nu}_m$, we get
\begin{eqnarray*}
\left(\V^{m+1}\right)^\top\MM_p\V^{m+1}&-&\left(\V^m\right)^\top \MM_p \V^m= \sum_{\nu=0}^{N_{sub}-1}\left(\left(\V^{\nu+1}_m\right)^\top\MM_p\V^{\nu+1}_m-\left(\V^{\nu}\right)^\top \MM_p \V^{\nu}  \right) \\
&\stackrel{\eqref{eq:disgrad_sub_v}}{=}& \sum_{\nu=0}^{N_{sub}-1} \left( (\V^{\nu+1}_m + \V^\nu_m)^\top
\MM_q  \int_{\tau_\nu}^{\tau_{\nu+1}}\LaB^1_i (\X_m(\tau))\du\tau \right) \left( \frac{\ce^m+\ce^{m+1}}{2} \right) \\
&\stackrel{\eqref{eq:dgsub_e}}{=}& -\left(\left(\ce^{m+1}\right)^\top \MM_{1} \ce^{m+1}-\left(\ce^{m}\right)^\top \MM_{1} \ce^{m}\right).
\end{eqnarray*}
Hence, the difference in the kinetic energy equals the negative difference in the potential energy (since the magnetic field is not changed in this system) and thus the total energy is conserved.

Moreover, Gauss' law is respected over time since we have in each subinterval $\tau \in [\tau_{\nu},\tau_{\nu}]$ that
\begin{equation*}
\X_m(\tau) = ((\tau_{\nu+1}-\tau) \X^\nu_m + (\tau-\tau_\nu)  \X^{\nu+1}_m)/\Delta \tau_\nu = \X^\nu_m +  \frac{\V_m^{\nu}+\V_m^{\nu+1}}{2} \left(\tau - \tau_{\nu}\right),
\end{equation*}
and hence 
\begin{equation*}
\frac{\du X_m(\tau)}{\du \tau} = \frac{\V_m^{\nu}+\V_m^{\nu+1}}{2},
\end{equation*}
which corresponds to the form of $\V(\tau)$ in \eqref{eq:dgsub_e}.

\begin{algorithm}
\caption{Second order conservative discrete gradient scheme with substepping.}\label{alg:vm_dgsub}
{\tiny
\For{$m=1, \ldots$}{
$\ce^{*} = S^{-1} \left( \left( \MM_{1}- \frac{\Delta t^2}{16} \C^T \MM_{2} \C\right) \ce^m + \Delta t \C^T \MM_{2} \cb^m \right)$ \;
$\cb^{*} = \cb^{m} - \frac{\Delta t}{4} \C \left(\ce^m+\ce^{*}\right)$\;
$\cj = 0$\;
$N = 0$\;
\For{$a=1, \ldots, N_p$}{
	$\vb_a^* = R(\bb^*) \vb_a^m$\;	
	$\xb_a^{\ii} = \xb_a^m + \frac{\Delta t}{2} \vb_a^*$\;
	$\vb_a^{**} = \vb^* + \frac{\Delta t}{2} \frac{q_s}{m_s} \LaB^1(\xb_a^{\ii}) \eb^*$\;
	$\cj = \cj + \Delta t q_s w_a \LaB^1(\xb_a^{\ii}) \vb_a$\;
	\For{$i=1,2,3$}{
	$N_i = N_i + \frac{\Delta t^2}{4} q_s w_a  \LaB_i^1(\xb_a^{\ii}) \LaB_i^1(\xb_a^{\ii})$ \;
	}
}
$\ce^{\ii} = \left(\MM_{1} - \begin{pmatrix}
N_1 & 0 & 0 \\ 0 & N_2 & 0 \\ 0 & 0 & N_3
\end{pmatrix}\right) \ce^* - \cj$\;
\For{$a=1, \ldots, N_p$}{
	$\vb_a^{\ii} = \vb_a^{**} + \frac{\Delta t}{2} \frac{q_s}{m_s} \LaB^1(\xb_a^{\ii}) \eb^{\ii}$\;
	$\xb_a^{\ii} = \xb_a^{\ii} + \frac{\Delta t}{2} \vb_a^{\ii}$\;
}
$residual=tolerance+1$\;
\While{$residual<tolerance$}{
$\cj = 0$\;
$\bar{\eb} = \frac{\eb^* + \eb^{\ii}}{2}$\;
\For{$a=1, \ldots, N_p$}{
	$\xb_a^{\text{before}} = \xb_a^*$\;
	$\vb_a^{\text{before}} = \vb_a^*$\;
	\For{$\nu=1,\ldots, N_{\nu}$}{
	$subresidual = subtolerance+1$\;
	$\xb_a^{\text{old}} = \xb_a^{\text{before}}$\;
	$\vb_a^{\text{old}} = \vb_a^{\text{before}}$\;
		\While{$subresidual < subtolerance$}{
			$\bar{\vb} =\frac{1}{2}\left( \vb_a^{\text{before}} + \vb^{\ii}\right)$\;
			$\xb_a^{\ii} = \xb_a^{\text{before}} + \Delta t \bar{\vb}$\;
			$\vb_a^{\ii} = \vb_a^{\text{before}} + \frac{\Delta t}{2} \int_{0}^1\frac{q_s}{m_s} \LaB^1(\xb_a^{\text{before}} + \tau \bar{\vb}) \du \tau \bar{\eb}$\;
			$subresidual = \max\left(\|\xb_a^{\ii}-\xb_a^{ \text{old}}\|_{\infty}, \|\vb_a^{\ii}-\vb_a^{ \text{old}}\|{\infty}\right)$\;
			
		$\xb_a^{\text{old}} = \xb_a^{\text{\ii}}$\;
		$\vb_a^{\text{old}} = \vb_a^{\text{\ii}}$\;
		}	
		$\cj =  \cj + \Delta t q_s w_a \int_0^1\LaB^1(\xb_a^{\text{old}} + \tau \bar{\vb}) \du \tau \bar{\vb}$\;
		$\xb_a^{\text{before}} = \xb_a^{\ii}$\;
		$\vb_a^{\text{before}} = \vb_a^{\ii}$\;
		}	
}
$\ce^{\text{old}} = \ce^{\ii}$\;
$\ce^{\ii} = \ce^* - M^{-1}\cj$\;
$residual = \|\ce^{\ii}-\ce^{\text{old}}\|_2$\;
}
\For{$a=1, \ldots, N_p$}{
	$\vb_a^{m+1} = R(\bb^*) \vb_a^{\ii}$\;
	$\xb_a^{m+1} = \xb_a^{\ii} + \frac{\Delta t}{2} \vb_a^{m+1}$\;
}
$\ce^{m+1} = S^{-1} \left( \left( \MM_{1}- \frac{\Delta t^2}{16} \C^T \MM_{2} \C\right) \ce^{\ii} + \Delta t \C^T \MM_{2} \cb^* \right)$ \;
$\cb^{m+1} = \cb^{*} - \frac{\Delta t}{4} \C \left(\ce^{\ii}+\ce^{m+1}\right)$\;
}
}
\end{algorithm}

\section{Numerical experiments}\label{sec:numerics}

In this section, we verify the conservation properties of our new time discretization methods for a number of test problems. We first study the reduced model in 1d2v phase-space with a perturbation along $x_1$, a magnetic field along $x_3$, and an electric field along the $x_1$ and $x_2$ directions. Moreover, we assume that the distribution function is independent of $v_3$. For this example, we report results on a Weibel instability, the two-stream instability, and a two-species simulation of an ion-accoustic wave. Also, we demonstrate the absence of the finite grid instablity. Finally, we show also results for the Weibel instability simulated in full 3d3v phase-space.

In all experiments reported, we have used a second oder Strang splitting method and third order splines for the 0-forms. For the iterative linear solvers, the tolerance is set to $10^{-15}$ and for the nonlinear iteration in the discrete gradient method we use a tolerance of $10^{-12}$ and a tolerance of $10^{-10}$ for the subiterations when they exist. Note that the tolerance of the linear solver is applied to the residual while the tolerance in the nonlinear iteration is directly applied to the fields. In order to balance the errors, we therefore use a more restrictive tolerance for the linear solvers. In the subiterations, only particle positions are involved which is why the tolerance can be chosen less restrictive. The implementation is based on the Fortran libraries SeLaLib\footnote{http://selalib.gforge.inria.fr/} and PLAF\footnote{http://jorek.gforge.inria.fr/documentation/plaf/html/index.html}.

\subsection{Weibel instability in 1d2v}


As a first example, we consider the Weibel instability \cite{Weibel:1959} in 1d2v. We use the same parameters that had already been considered with the Hamiltonian splitting time discretization in \cite{Kraus17}. The initial distribution and fields are of the form
\begin{align*}
\f(x, \vb,t=0) &= \frac{1}{2\pi \sigma_1 \sigma_2} \exp \left(- \frac{1}{2}\left( \frac{v_1^2}{\sigma_1^2} + \frac{v_2^2}{\sigma_2^2} \right) \right) \left( 1 + \alpha \cos( kx)\right), \quad x \in [0,2\pi/k),\\
B_3(x,t=0) &= \beta \cos(kx),\\
E_2(x,t=0) &= 0,
\end{align*}
and $E_1(x,t=0)$ is computed from Poisson's equation. 
In our simulations, we use the following choice of parameters, $\sigma_1 = 0.02/\sqrt{2}$, $\sigma_2 = \sqrt{12} \sigma_1$, $k=1.25$, $\alpha = 0$ and $\beta = 10^{-4}$ and simulate for 500 time units. We use 100,000 particles and 32 grid points.

 We run the simulation with the average-vector-field method and the Gauss-conserving discrete gradient method for various time steps and compare to the Hamiltonian splitting. Table \ref{tab:weibel_1d2v} shows the conservation properties of the various runs. The numerical experiments verify energy conservation of the new semi-implicit methods and conservation of Gauss' law for the discrete gradient method and the Hamiltonian splitting. In particular, we note that the error in Gauss' law is satisfied up to round-off errors in the Hamiltonian splitting as well as the conservative discrete gradient method. Moreover, the implicit methods conserve energy to the accuracy of the field solver.  We can also see that the semi-implicit methods allow for larger time steps where the Hamiltonian splitting becomes unstable due to the stability constraint $\Delta t \leq \sqrt{\frac{17}{42}}\Delta x \approx 0.099935$ for the explicit scheme for Maxwell's equations.
 
  The number of nonlinear iterations in the discrete gradient scheme are about 4 ($\Delta t=0.025$), 5 ($\Delta t=0.05$), 6 ($\Delta t=0.1$), and 8 ($\Delta t=0.2$). Hence, we see only a moderate increase in the number of iterations needed. Roughly speaking the cost of the implicit method is a factor ``number of iterations'' more expensive than the explicit Hamiltonian splitting. This shows that the computational costs of our implicit and semi-implicit schemes are comparable to the costs that are required for the EC-PIC and EC-SIM algorithms (cf.~\cite{Lapenta17}).

 In this example, the error of the various methods is comparable as examplified for the maximum error over time in the magnetic energy shown in the last column of Table \ref{tab:weibel_1d2v}. The error in the magnetic energy is computed compared to a solution with a time step $\Delta t = 0.0125$ and the same method.

\begin{table}
\caption{Weibel instability in 1d2v phase space: Comparison of the conservation properties for various integrators. For the average-vector-field and the discrete gradient method, we show results for two different oderings of the individual operators (the first set of experiments refers to the case with the ordering shown in the Algorithms).} \label{tab:weibel_1d2v}
\centering
\begin{tabular}{|c|c|c|c| c |}
\hline
method & $\Delta t$ & Gauss & energy & error magn. energy\\
\hline \hline
HS & $0.025$ & 2.39E-15 & 3.48E-07 & 5.12E-06 \\
HS & $0.05$  & 2.39E-15 & 1.39E-06 & 1.88E-05 \\
HS & $0.1$  &  --- & --- & --- \\
\hline 
AVF (O3,O1,O2,O4) & 0.025 & 1.25E-07 & 6.599E-15 & 1.45E-05 \\
AVF (O3,O1,O2,O4)  & 0.05 & 2.68E-07 & 3.08E-14 & 1.53E-05 \\
AVF (O3,O1,O2,O4)  & 0.1 & 1.10E-06 & 1.11E-14 & 4.15E-05 \\
AVF (O3,O1,O2,O4)  & 0.2 & 5.80E-06 & 8.76E-15 & 1.24E-04 \\
\hline 
AVF (O1,O2,O3,O4) & 0.025 & 8.77E-08 & 7.44E-15 & 5.06E-06 \\
AVF (O1,O2,O3,O4)  & 0.05 & 3.38E-07 & 1.97E-14 & 2.02E-05 \\
AVF (O1,O2,O3,O4)  & 0.1 & 1.93E-06 & 8.81E-15 & 4.06E-05 \\
AVF (O1,O2,O3,O4)  & 0.2 & 7.23E-06 & 1.03E-14 & 1.12E-04 \\
\hline 
DisGrad (O3,O2,O4) & 0.025 & 2.28E-15 & 1.56E-13 & 8.22E-06\\
DisGrad (O3,O2,O4)  & 0.05 & 2.32E-15 & 2.60E-14 & 1.50E-05\\ 
DisGrad  (O3,O2,O4) & 0.1 & 2.14E-15 & 1.66E-13 & 4.11E-05 \\
DisGrad  (O3,O2,O4) & 0.2 & 2.09E-15 & 5.48E-14 & 1.20E-04\\
\hline
DisGrad (O2,O3,O4) & 0.025 & 2.72E-15 & 1.01E-13  & 5.68E-06 \\
DisGrad (O2,O3,O4)  & 0.05 & 2.23E-15 &  5.92E-15 & 1.69E-05\\
DisGrad (O2,O3,O4)  & 0.1 & 2.24E-15 &  1.81E-13 & 3.99E-05\\
DisGrad (O2,O3,O4)  & 0.2 & 2.24E-15 & 2.15E-14 & 1.09E-04  \\
\hline
\end{tabular}
\end{table}

\subsection{Two-stream instability}


As a second example, we look at a classical electrostatic test case known as the two-stream instability with the following initial value
\begin{align*}
\f(x, \vb,t=0) &= \frac{1}{4\pi } \left( \e^{ - 0.5(v_1-2.4)^2} +  \e^{ - 0.5(v_1+2.4}\right) 
 \e^{ - 0.5v_2^2}, \quad x \in [0,10\pi),\\
B_3(x,t=0) &= 0,\\
E_2(x,t=0) &= 0,
\end{align*}
and $E_1$ computed from Poisson's equation. We have simulated the problem over 200 time units with 64 grid points and 64,000 particles with the various integrators for time steps of 0.025, 0.05, 0.1, 0.2, and 0.4. Again the time step for the explicit scheme is restricted due to the Maxwell's equations by $\Delta t \leq \sqrt{\frac{17}{42}} \Delta x \approx 0.19882$. We note that larger time steps would therefore be possible for this electrostatic example if we solve the Vlasov--Amp{\`e}re equation instead. For the average vector field and the discrete gradient method the operators are ordered as shown in the algorithms. The simulations with the discrete gradient method took the following number of iterations on average: 7 ($\Delta t = 0.025$), 7 ($\Delta t = 0.05$), 9 ($\Delta t = 0.1$), 11 ($\Delta t = 0.1$), 16 ($\Delta t = 0.4$).

\begin{figure}
\centering
\includegraphics[scale=1.0]{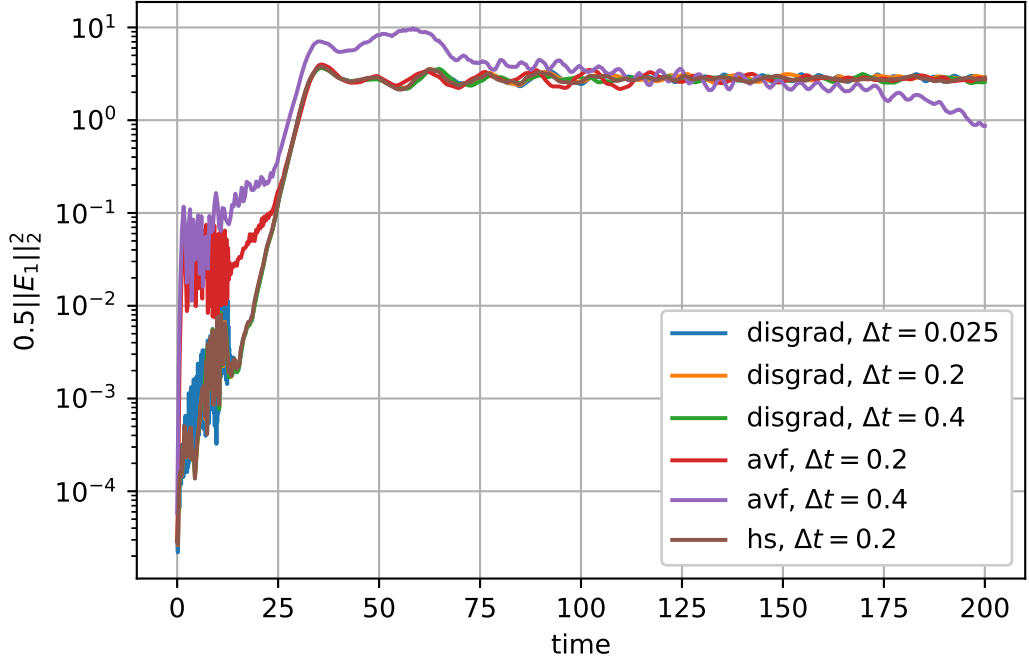}
\caption{Two-stream instability: Time evolution of the first component of the electric energy for various configurations.}\label{fig:tsi}
\end{figure}

The conservation properties of the various simulations are summarized in Table \ref{tab:tsi}. All simulations with time step smaller than 0.4 capture the linear growth rate quite accurately and show only small variations in the nonlinear phase. For a time step of 0.4, the stability condition of the Hamiltonian splitting is violated, the average vector field was still able to capture the linear growth rate but showed large deviation in the nonlinear phase while the results of the discrete gradient scheme remains very accurate. Moreover, the initially random electric field deviates considerably for the average-vector field method with time steps of 0.2 and 0.4 which can be explained by the fact that Gauss' law is violated to a very large degree in these simulations. Figure \ref{fig:tsi} shows the electric energy as a function of time for the discrete gradient method at a time step of 0.025 as a reference and the simulations with $\Delta t=0.2$ and $\Delta t = 0.4$. 

Even though we have demonstrated that large time steps are possible with the implicit methods for the examples of the two-stream and the Weibel instability, large time steps are not particularly beneficial for these examples due to accuracy. As a next example we will therefore study a multi-species problem where physical effects take place on the time-scale of the ion so that large time steps are interesting for the electron dynamics.

\begin{table}
\caption{Two-stream instability 1d2v phase space: Comparison of the conservation properties for various integrators.} \label{tab:tsi}
\centering
\begin{tabular}{|c|c|c|c| }
\hline
method & $\Delta t$ & Gauss & energy \\
\hline \hline
HS & $0.025$ &4.77E-15 &7.37E-04 \\
HS & $0.05$  &4.11E-15 & 2.79E-03  \\
HS & $0.1$  &  4.77E-15 & 9.37E-03\\
HS & $0.2$ & 4.44E-15 & 4.81E-01 \\
HS & $0.4$ & --- & ---\\
\hline 
AVF & $0.025$ & 1.26E-03 & 6.71E-12\\
AVF & $0.05$  & 5.17E-03 & 5.68E-12 \\
AVF & $0.1$  & 2.56E-02 & 6.08E-12 \\
AVF & $0.2$ & 2.07E-01 & 5.49E-12 \\
AVF & $0.4$ & 1.30E-00 & 5.68E-12\\
\hline 

DisGrad & 0.025 & 4.55E-15 &  9.69E-12\\
DisGrad & 0.05 &  3.89E-15 & 8.38E-11 \\
DisGrad & 0.1 & 3.55E-15 & 5.68E-12 \\
DisGrad & 0.2 & 4.33E-15 & 1.68E-11 \\
DisGrad & 0.4 &  4.72E-15 & 2.71E-11\\
\hline
\end{tabular}
\end{table}

\subsection{Ion acoustic wave}


As a third example, we consider the ion acoustic wave excited by an ion density perturbation. The example is electrostatic and involves electrons and ions. We normalize mass and temperature to their values relative to the electrons. The initial distributions are then given by
\begin{equation*}\begin{aligned}
f_e(x,v,t=0) &= \frac{1}{\sqrt{2\pi}} \exp \left( -\frac{v^2}{2} \right), \\
f_i(x,v,t=0) &= \frac{1}{\sqrt{2\pi \frac{T_i}{m_i}}} \exp \left( -\frac{v^2}{2 \sqrt{\frac{T_i}{m_i}}} \right) \left( 1 + \alpha \cos \left( \frac{2 \pi}{L}x\right) \right). 
\end{aligned}\end{equation*}
We use the following parameters $T_i = 10^{-4}$, $m_i = 200$, $\alpha= 0.2$, $L=10$.  

The ion acoustic wave problem is solved with 32 grid points and 128,000 particles per species with the explicit Hamiltonian splitting, the average vector field method, and the discrete gradient method with and without substepping for various time steps. In this case, we use the ordering that places operator 3 after operators 1 and 2 since this gives more accurate results in the present example. Table \ref{tab:iaw_cons} shows the conservation properties of the various methods for different time steps. As expected Gauss' law is satisfied to machine precision for the Hamiltonian splitting as well as the discrete gradient method. On the other hand, the error in Gauss' law is considerable for the average-vector-field method and increases with the time step. The average-vector-field method and the discrete gradient method conserve energy---in contrast to the Hamilitonian splitting---up to the tolerance of the iterative solvers.  

\begin{figure}
\centering
\includegraphics[scale=0.5]{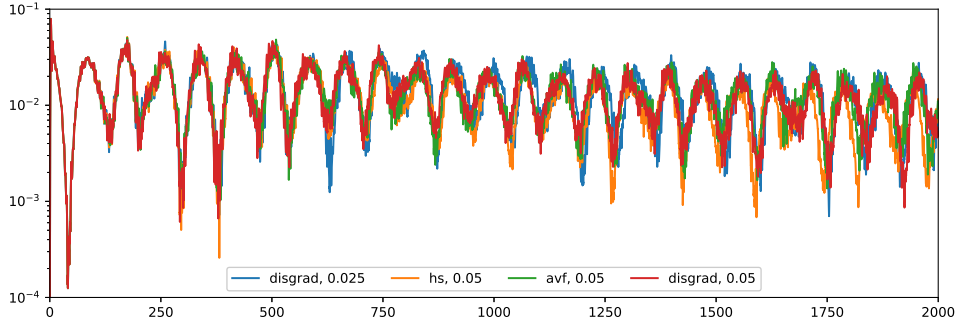}
\caption{Ion acoustic wave: Time evolution of the first component of the electric energy with various time propagation schemes at $\Delta t = 0.05$ compared to a reference simulation with $\Delta t= 0.025$.}\label{fig:iaw_dt005}
\end{figure}

\begin{figure}
\centering
\includegraphics[scale=0.5]{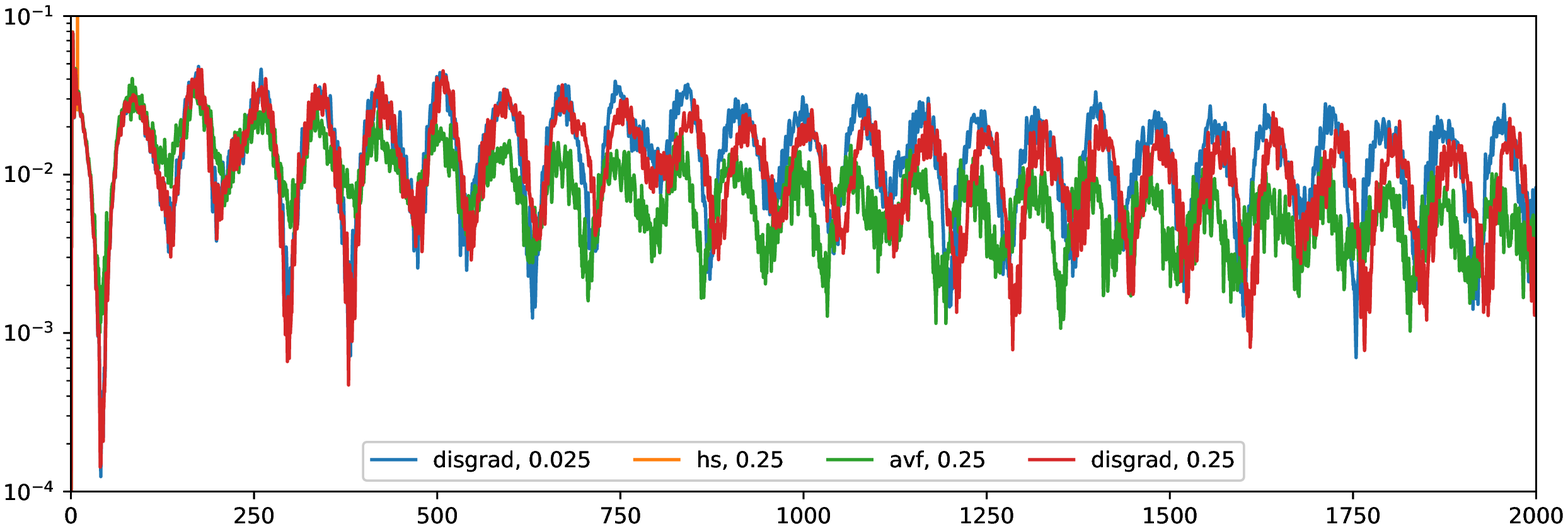}
\caption{Ion acoustic wave: Time evolution of the first component of the electric energy with various time propagation schemes at $\Delta t = 0.25$ compared to a reference simulation with $\Delta t= 0.025$.}\label{fig:iaw_dt025}
\end{figure}

\begin{figure}
\centering
\includegraphics[scale=0.5]{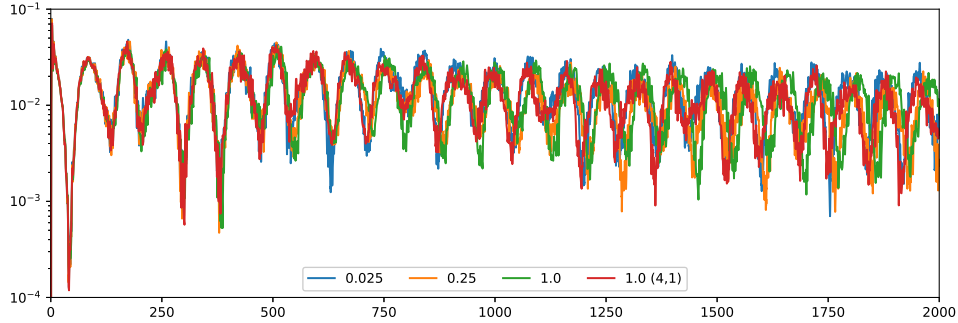}
\caption{Ion acoustic wave: Time evolution of the first component of the electric energy comparing the discrete gradient method without sub-stepping and a time step of $\Delta t=0.025$, $\Delta t=0.25$, and $\Delta t=1.0$ compared to the discrete gradient method with $\Delta t=1.0$ and four substeps for the electrons.}\label{fig:iaw_dg}
\end{figure}

In order to judge the quality of the solution, we look at the evolution of the first component of the electric energy over time. As a reference, we show the solution of the discrete gradient method with a time step of $\Delta t=0.025$. Figure \ref{fig:iaw_dt005} shows the results for the Hamiltonian splitting, the average-vector-field, and the discrete gradient solution with a time step of $\Delta t = 0.05$. We can see that all three time stepping schemes give quite good results for this relatively small time step. Next we increase the time step to $\Delta t = 0.25$ in Figure \ref{fig:iaw_dt025} where the explicit method become unstable and the discrete gradient method gives clearly better results than the average-vector field method. Finally, Figure~\ref{fig:iaw_dg} shows that the solution of the discrete gradient method becomes clearly worse for $\Delta t = 1.0$. When introducing a substepping scheme with four substeps for the electrons only, the quality of the solution at $\Delta t=1$  and $\Delta \tau = 0.25$ for the electrons is comparable to the discrete gradient method with a total time step of $\Delta t = 0.25$ for both species.

The efficiency of the iterative discrete gradient method is hampered by the number of iterations needed for the nonlinear solution to converge. In particular, for the smallest time step considered ($\Delta t = 0.025$), we already need about 6 iterations on average which renders the method uncompetitive compared to the explicit one. When increasing the time step, the iteration count increases to about 7 iterations for $\Delta t=0.05$, 12 iterations for $\Delta t =0.25$ , 17 for $\Delta t = 0.5$, and 32 for $\Delta t = 1.0$. We can see that the increase of iterations needed between $\Delta t = 0.05$ and $\Delta t = 0.25$ is quite small compared to the increase in the time step. Then, the iteration count starts increasing at a higher rate. However, the quality of the solution with such large time steps is not very good either so that the high value of the iteration tolerance is questionable in those simulations. The substepping method still produces good results for $\Delta t = 1.0$ and four substeps for the electrons and no substeps for the ions. In this case, about 16 outer iterations are needed per time step and about 5 inner iterations for both electrons (per substep) and ions.
\begin{table}
\caption{Ion acoustic wave: Comparison of the conservation properties for various integrators.} \label{tab:iaw_cons}
\centering
\begin{tabular}{|c|c|c|c|}
\hline
method & $\Delta t$ & Gauss & energy\\
\hline \hline
HS & 0.025 &  9.10E-15 & 1.57E-05 \\
HS & 0.05 & 9.96E-15 & 8.93E-05 \\
\hline
AVF & 0.025 &   3.07E-4 & 4.28E-13\\
AVF & 0.05 & 1.78E-03 & 3.27E-13 \\
AVF & 0.25 & 7.78E-02 & 2.75E-13\\
AVF & 0.5 & 2.31E-01 & 2.58E-13\\
AVF & 1.0 & 3.93E-01 & 2.36E-13 \\
\hline
DisGrad & 0.025 & 1.11E-14 & 9.34E-12 \\
DisGrad & 0.05 & 1.22E-14 & 2.96E-12 \\
DisGrad & 0.25 & 1.28E-14 & 2.96E-12 \\
DisGrad & 0.5 & 1.49E-14 & 4.52E-12 \\
DisGrad & 1.0 & 1.53E-14 & 1.85E-11 \\
\hline
DisGrad, sub(4,1) & 1.0 &  1.71E-14 & 5.45E-13\\
\hline
\end{tabular}
\end{table}

Note that a linear dispersion analysis shows that the electric energy should oscillate with no damping in this case. Our results, however, show a slight damping. This is not an effect of the time stepping scheme but of the spatial discretization. Adding more particles or smoothing the fields improves the results.

\subsection{Finite grid instability}

\begin{figure}
\centering
\includegraphics[scale=1.0]{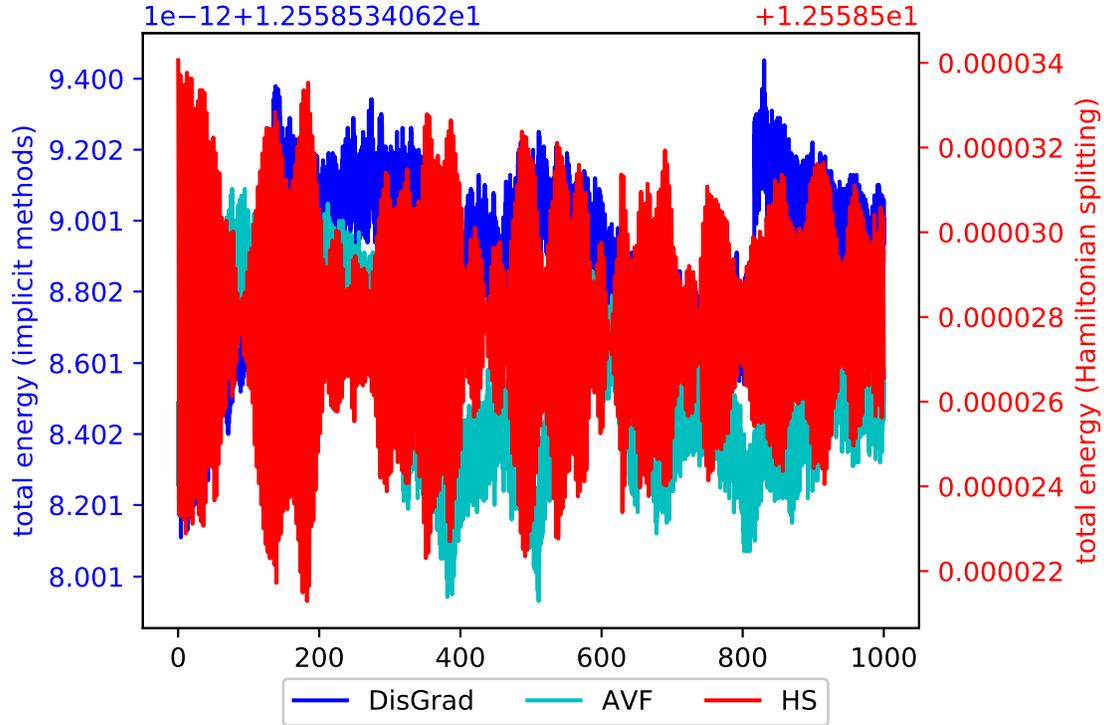}
\caption{Finite grid instability: Time evolution of the total energy for simulation with resolution $>10$ times the Debye length.}\label{fig:fgi}
\end{figure}

For numerical schemes that lack energy conservation, artificial heating occurs if the Debye length is not resolved on the spatial grid. Let us consider a Maxwellian initial condition of the form
\begin{equation*}
f(x,v_1,v_2) = \frac{1}{2\pi \sigma} \exp\left( - 0.5 \left(\frac{v_1^2+v_2^2}{\sigma^2} \right)\right)
\end{equation*}
in a periodic box of length $L=50 \pi$. In our normalized units the plasma frequency is equal to one and the Debye length then takes the value $\sigma$ of the thermal velocity. We choose a spatial resolution of 64 grid points, i.e. $\Delta x \approx 2.4$, and a thermal velocity of $\sigma = 0.2$. In this case, the Debye length is a factor 10 larger than the spatial resolution and a standard explicit method would suffer from the finite grid instability. On the contrary, we do not see a grid heating for any of our methods. Figure \ref{fig:fgi} shows the evolution of the total energy over time for simulations with a very small time step of $\Delta t = 0.05$. In particular, the energy is conserved over time up to the solver tolerance of $10^{-12}$ for the implicit method. For the explicit Hamiltonian splitting the energy is not conserved but it shows an oscillatory behavior instead of the steady increase that is referred to as finite grid instability. We note that the energy error of the Hamiltonian splitting is solely caused by the splitting in time and its error is a function of $\Delta t$ only. The spatial discretization is energy-conserving which is why this explicit scheme does not suffer from the finite grid instability.

\subsection{Weibel instability in 3d3v}


In order to demonstrate the potential of the algorithm also in full phase space, we consider the Weibel test case in a simulation in the full six-dimensional phase space. The initial distribution is of the form
\begin{align*}
\f(\xb, \vb,t=0) &= \frac{1}{(2\pi)^{3/2} \sigma_1 \sigma_2^2} \exp \left(- \frac{1}{2}\left( \frac{v_1^2}{\sigma_1^2} + \frac{v_2^2+v_3^2}{\sigma_2^2} \right) \right) , \quad \xb \in [0,2\pi/k)^3,\\
\Bb(\xb,t=0) &= ( 0,\, 0,\, \beta \cos(kx_1))^\top,
\end{align*}
and the electric field at time zero is computed from Poisson's equation. In our simulation, we choose the parameters $\sigma_1 = 0.02/\sqrt{2}$, $\sigma_2 = \sqrt{12} \sigma_1$, $k=1.25$, and  $\beta = 0.01$ and simulate for 250 time units. Note that the problem is mainly depending on the variable $x_1$ which is why we resolve the spatial grid more along this direction, namely we choose a resolution of $16 \times 8 \times 8$ grid points. The simulation uses $N_p=100,000$ particles. 

Table \ref{tab:weibel_3d3v} shows the conservation properties of the various integrators with different time steps. The results here are shown for the standard ordering of the operators that maximizes the arithmetic intensity. Again the conservation properties are verified. For this case, the geometric discrete gradient method needs about 5 ($\Delta t = 0.05$), 6 ($\Delta t = 0.1$), 9 ($\Delta t = 0.2$), or 14 ($\Delta t = 0.4$) iterations per time step on average. The error in the magnetic energy is computed compared to a solution with the same method and a time step of $\Delta t=0.0125$. The accuracy is comparable for all methods with a slight disadvantage for the Hamiltonian splitting.

\begin{table}
\caption{Weibel instability in 3d3v phase space: Comparison of the conservation properties for various integrators.} \label{tab:weibel_3d3v}
\centering
\begin{tabular}{|c|c|c|c| c |}
\hline
method & $\Delta t$ & Gauss & energy & error magn. energy\\
\hline \hline
HS & $0.05$ & 5.80E-16 & 6.70E-06 & 1.23E-03 \\
HS & $0.1$  & 5.70E-16 & 2.72E-05 & 5.03E-03 \\
HS & $0.2$  &  --- & --- & --- \\
\hline 
AVF & 0.05 & 3.47E-07 & 2.59E-13 & 8674E-04 \\
AVF & 0.1 & 1.60E-06 & 3.76E-14 & 3.63E-03 \\
AVF & 0.2 & 6.00E-06 & 2.85E-14 & 1.35E-02 \\
AVF & 0.4 & 2.27E-05 & 6.88E-14 & 2.91E-02 \\
\hline 
DisGrad & 0.05 & 1.27E-15 & 3.53E-13 &8.61E-04\\
DisGrad & 0.1 & 1.20E-15 & 1.50E-14 & 3.63E-03\\
DisGrad & 0.2 & 1.06E-15 & 2.22E-13 & 1.33E-02 \\
DisGrad & 0.4 & 1.03E-15 & 8.22E-15 & 2.92E-02\\
\hline
\end{tabular}
\end{table}


\section{Conclusions and outlook}

We have described a general procedure to derive energy-conserving time-stepping methods for the geometric electromagnetic particle-in-cell method based on discrete gradients and an antisymmetric splitting of the Poisson matrix. In particular, we derived a semi-implicit scheme based on the average-vector-field method. The method yields good results with little computational overhead per time step compared to the Hamiltonian splitting when conservation of Gauss' law is not critical and, in the same time, allows for larger time steps. Moreover, we have derived an implicit method that conserves both energy and Gauss' law. Due to the fact that a nonlinear iteration that couples the particle and field degrees of freedom is necessary the method comes with a considerable computational overhead. On the other hand, combined with a substepping for the fast species, it yields quite accurate results for large time steps and has therefore the potential to be more efficient in multiscale simulations with realistic mass ratio between electron and ion species. Also the proposed (semi-)implicit methods have the  potential to be extended to more complex geometries where splitting of the components of the $H_p$ in the Hamiltonian splitting is not possible anymore. Finally, we have also shown the absense of the finite grid instability for the GEMPIC semi-discretization independent of the time stepping scheme.

%
%
%

\section*{Acknowledgement} The authors acknowledge discussions with Michael Kraus, Philip J. Morrison, Benedikt Perse, Jalal Lakhlili, and Ahmed Ratnani. This work has been carried out within the framework of the EUROfusion Consortium and has received funding from the Euratom research and training program 2014-2018 and 2019-2020 under grant agreement No 633053. The views and opinions expressed herein do not necessarily reflect those of the European Commision.

\appendix

\section{Stability analysis}

A comprehensive stability analysis of the geometric particle-in-cell method is quite hard to achieve due to the fact that the scheme is highly nonlinear. A rather general stability analysis for electromagnetic particle-in-cell schemes was provided by Godfrey \cite{Godfrey75}, however, it does not directly apply to our methods since it is based on a staggered time step while we use a more complex splitting of the equations. On the other hand, the stability limits of explicit particle-in-cell schemes are mostly related to the way the curl-part in Maxwell's equations is solved and the propagation of electrostatic Langmuir waves. In this section, we will therefore perform a stability analysis for the curl-part in Maxwell's equation in 1d and for Langmuir waves. 

\subsection{Stability analysis for Langmuir waves}

In this section, we perform a von Neumann stability analysis for Langmuir waves as in \cite{Lapenta17}. For the electrostatic case, both implicit schemes in a Strang splitting combination read as follows
\begin{eqnarray}\label{avf:prop_es}
	\X^{m+1/2} = \X^{m-1/2} + \Delta t \V^m, \quad \V^{m+1} = \V^m + \frac{\Delta t}{2}\frac{q_s}{m_s}\left( \E^m + \E^{m+1} \right). 
\end{eqnarray}
For the von Neumann stability analysis, the time evolution is assumed to be harmonic, i.e.
\begin{eqnarray*}
\xb_p^{m} = \bar \xb_p \e^{\im \omega m/\Delta t}, \quad \vb_p^{m} = \bar \vb_p \e^{\im \omega m/\Delta t},\quad \Eb_p^{m} = \bar \Eb_p \e^{\im \omega m/\Delta t}.
\end{eqnarray*}
Inserting the time harmonic ansatz into the scheme \eqref{avf:prop_es}, we find
\begin{eqnarray*}
&&\bar \xb_p 2\im \sin\left(\frac{\Delta t}{2}\omega\right) = \Delta t \bar \vb_p, \\
&&\bar \vb_p 2\im \sin\left(\frac{\Delta t}{2}\omega\right) = \frac{q_s}{m_s} \Delta t \cos\left(\frac{\Delta t}{2}\omega\right).
\end{eqnarray*}
In Fourier space the electric field for cold plasma Langmuir waves can be related to the displacement by the electron plasma frequency $\omega_{pe}$ as (cf. \cite{Lapenta17})
\begin{equation*}
\frac{q_s}{m_s} \bar \Eb = - C^2\omega_{pe}^2 \bar \xb,
\end{equation*}
for some constant $C$. With this assumption, we obtain
\begin{eqnarray*}
&&\bar \xb_p 2\im \sin\left(\frac{\Delta t}{2}\omega\right) - \Delta t \bar \vb_p = 0, \\
&&  \Delta t \cos\left(\frac{\Delta t}{2}\omega\right) \omega_{pe}^2 C^2  \xb_p +\bar \vb_p 2\im \sin\left(\frac{\Delta t}{2}\omega\right) = 0.
\end{eqnarray*}
This system of linear equations has a solution if the determinant of the matrix vanishes, i.e. if
\begin{equation*}
	-4 \sin^2\left(\frac{\Delta t}{2}\omega\right) + C^2\Delta t^2 \omega_{pe}^2\cos\left(\frac{\Delta t}{2}\omega\right) = 0.
\end{equation*}
For this case, we hence get the same stability condition as for the semi-implicit scheme introduced by Lapenta in \cite{Lapenta17}. The equation has real solutions independent of $\Delta t$ and, hence, the scheme is unconditionally stable.

The Hamiltonian splitting method with second order Strang splitting, on the other hand, yields the standard explicit leap frog scheme for electrostatics
\begin{equation*}
	\X^{m+1/2} = \X^{m-1/2} + \Delta t \V^m, \quad \V^{m+1} = \V^m + \Delta t\frac{q_s}{m_s}\Eb^{m+1/2}. 
\end{equation*}
In this case the determinant condition reads
\begin{equation*}
	-4 \sin^2\left(\frac{\Delta t}{2}\omega\right) + C^2\Delta t^2 \omega_{pe}^2 = 0.
\end{equation*}
For $C\omega_{pe} \Delta t> 2$, there are only complex conjugate solutions, i.e. we have one growing solution which leads to numerical instabilities.

\subsection{Curl-part of Maxwell's equations}

Let us now consider the curl-part of Maxwell's equations in 1d with a second order Strang splitting propagator. Let us denote $\tilde{\eb} = \MM_1 \eb$ and $\tilde{\bb} = \MM_2 \bb$. With this notation the explicit version of the Maxwell equations with a finite element description reads for degree of freedom $j$
\begin{equation*}\begin{aligned}
\me_j^{n+1/2} &=& \me_j^n + \frac{\alpha}{2} \left( \mb_j^{n}-\mb_{j+1}^{n}\right),\\
b_j^{n+1} &=& b_j^{n} - \alpha \left( e_j^{n+1/2}-e_{j-1}^{n+1/2}\right), \\
\me_j^{n+1} &=& \me_j^{n+1/2} + \frac{\alpha}{2} \left( \mb_j^{n+1}-\mb_{j+1}^{n+1}\right).
\end{aligned}\end{equation*}
where $\alpha = \frac{\Delta t}{\Delta x}$.

For one Fourier mode $k$, we use the ansatz 
\begin{equation*}
e_j^n = \bar{e} \xi^n \exp(\im k x_j), b_j^n = \bar{b} \xi^n \exp(\im k x_j).
\end{equation*}
Then, we have the following relation after multiplication with the mass matrices
\begin{equation*}
\me_j^n = \lambda^{(p)}_k e_j^n, \quad \mb_j^n = \lambda^{(p-1)}_k b_j^n,
\end{equation*}
where $\lambda^{(q)}_k$ denotes the $k$th eigenvalue of the mass matrix for $q$th order splines. The Fourier transformed difference equations then have the following form
\begin{equation*}
\underbrace{\begin{pmatrix}
(\xi-1)\lambda^{(p)}_k & -(\xi+1)\frac{\alpha}{2} \left(1-\e^{\im \Delta x k} \right) \lambda^{(p-1)}_k\\
\alpha \left(1-\e^{-\im \Delta x k} \right)\lambda^{(p)}_k & (\xi-1) \lambda^{(p)}_k + \alpha^2  \lambda^{(p-1)}_k \left(1-\cos( \Delta x k)\right)
\end{pmatrix}}_{:=D} \begin{pmatrix}
\bar{e} \\ \bar{b}
\end{pmatrix} = 0.
\end{equation*}
To find a solution, we compute $\xi$ such that the determinant is zero
\begin{equation*}
\det(D) = \lambda^{(p)}_k\left(\lambda^{(p)}_k (\xi-1)^2 + \alpha^2 2(1-\cos(k\Delta x))\xi \lambda^{(p-1)}_k \right)= \lambda^{(p)}_k\left(\xi^2 -2q\xi+1\right),
\end{equation*}
where $q = 1-\frac{\lambda^{(p-1)}_k}{\lambda^{(p)}_k}\alpha^2 \left(1-\cos(k\Delta x)\right)$.
The roots of the equation $\det(D)=0$ can be expressed as
\begin{equation*}
\xi_{+/-} = q \pm \sqrt{q^2-1}.
\end{equation*}
For stability, we need to have $|\xi|\leq 1$ and thus $|q| \leq 1$ which yields the condition
\begin{equation*}
0 \leq \frac{\lambda^{(p-1)}_k}{\lambda^{(p)}_k}\alpha^2 \left(1-\cos(k\Delta x)\right) \leq 2,
\end{equation*}
for all values of $k$. For $p=1,2,3$ this yields the following conditions on $\alpha$:
\begin{itemize}
\item $p=1$: $\alpha \leq \sqrt{\frac{1}{3}}$,
\item $p=2$: $\alpha \leq \sqrt{\frac{2}{5}}$,
\item $p=3$: $\alpha \leq \sqrt{\frac{17}{42}}$.
\end{itemize}
We note that a Lie splitting would yield the same determinant (up to the multiplicative factor $\lambda^{(p)}_k$) and, hence, yields the same stability limit.

Next, we consider the implicit variant. Since we will show that the scheme in unconditionally stable, it suffices to consider the Lie splitting
\begin{equation*}
\begin{aligned}
b_j^{n+1} = b_j^{n} - \frac{\alpha}{2} \left( e_j^n+e_j^{n+1}-e_{j-1}^n-e_{j-1}^{n+1}\right), \\
\me_j^{n+1} = \me_j^n + \frac{\alpha}{2} \left( \mb_j^{n}+\mb_j^{n+1}-\mb_{j-1}^{n}-\mb_{j-1}^{n+1}\right).
\end{aligned}\end{equation*}
With the same ansatz, we now get the following equation for mode $k$
\begin{equation*}
\begin{pmatrix}
\xi-1 & \frac{\alpha}{2} \left(1-\e^{-ik\Delta x} \right)(\xi+1) \\
- \frac{\alpha}{2}\left(1-\e^{ik\Delta x} \right)\lambda^{(p-1)}_k(\xi+1) & (\xi-1) \lambda^{(p)}_k
\end{pmatrix} \begin{pmatrix}
\bar{b} \\ \bar{e}
\end{pmatrix} = 0.
\end{equation*}
This yields the following expression for the determinant
\begin{equation*}\begin{aligned}
\det (D) =& \left(\lambda^{(p)}_k + \frac{\alpha^2}{2} (1-\cos(k\Delta x)) \lambda^{(p-1)}_k \right) \xi^2 + \left(\lambda^{(p)}_k - \frac{\alpha^2}{2} (1-\cos(k\Delta x)) \lambda^{(p-1)}_k \right) \xi +\\
& \left(\lambda^{(p)}_k + \frac{\alpha^2}{2} (1-\cos(k\Delta x)) \lambda^{(p-1)}_k \right).
\end{aligned}\end{equation*}
Solving the equation $\det (D) = 0$ for $\xi$ yields 
\begin{equation*}
\xi_{+/-} = q \pm \sqrt{q^2-1},
\end{equation*}
where $q = \frac{\lambda^{(p)}_k - \frac{\alpha^2}{2} (1-\cos(k\Delta x)) \lambda^{(p-1)}_k}{\lambda^{(p)}_k + \frac{\alpha^2}{2} (1-\cos(k\Delta x)) \lambda^{(p-1)}_k}$. In this case, it holds that $q^2-1 \leq 0$ independent of $k$, $\Delta x$, $p$ and $\alpha$ since both $\lambda^{(p)}_k\geq 0$ and  $\frac{\alpha^2}{2} (1-\cos(k\Delta x)) \lambda^{(p-1)}_k \geq 0$. Hence, the scheme is unconditionally stable.

\end{document}